\documentclass{article}
\usepackage[utf8]{inputenc}
\usepackage{graphicx}
\usepackage{amsmath}
\usepackage{amsfonts}
\usepackage{mathrsfs}
\usepackage{mathtools}
\usepackage[top=45truemm,bottom=45truemm,left=30truemm,right=30truemm]{geometry}
\usepackage{here}
\usepackage{amsthm}
\usepackage{algorithmic}
\usepackage{algorithm}
\usepackage{url}
\usepackage{stmaryrd}
\usepackage[title]{appendix}
\usepackage{amssymb}
\usepackage[all]{xy}
\usepackage{amscd}
\usepackage{color}
\usepackage{extarrows}
\usepackage{multirow}
\usepackage{enumitem}
\usepackage[
        dvipdfmx,
        hypertexnames=false,
        hyperindex,
        bookmarks=false,
        colorlinks=true,
        linkcolor=magenta,
        citecolor=magenta,
        urlcolor=blue
]{hyperref}
\usepackage{bm}

\newtheorem{theorem}{Theorem}[subsection]
\newtheorem{lemma}[theorem]{Lemma}
\newtheorem{proposition}[theorem]{Proposition}
\newtheorem{corollary}[theorem]{Corollary}
\theoremstyle{definition}
\newtheorem{definition}[theorem]{Definition}
\newtheorem{example}[theorem]{Example}
\newtheorem{remark}[theorem]{Remark}
\newtheorem{condition}[theorem]{Condition}

\theoremstyle{plain}
\newtheorem{theor}{Theorem}
\theoremstyle{plain}

\theoremstyle{plain}

\theoremstyle{plain}

\numberwithin{equation}{subsection}

\title{Degree bounds and synchronization in Gr\"{o}bner basis computations for affine semi-regular systems}
\date{\today}
\author{Momonari Kudo\thanks{Fukuoka Institute of Technology, 3-30-1 Wajiro-higashi, Higashi-ku, Fukuoka, 811-0295 Japan, \texttt{m-kudo@fit.ac.jp}} \and Kazuhiro Yokoyama\thanks{Rikkyo University, 3-34-1 Nishi-Ikebukuro, Toshima-ku, Tokyo, 171-8501 Japan, \texttt{kazuhiro@rikkyo.ac.jp}}}
\begin{document}

\maketitle

\abstract{
Determining the complexity of computing Gr\"obner bases is an important problem in both theory and practice, and solving degrees provide a central measure of this complexity.
We study solving degrees and Gr\"obner basis computations for affine polynomial systems, with particular emphasis on semi-regular sequences.

We first derive two upper bounds for the maximum Gr\"obner basis degree of the homogenized system.
One is based on a regular initial subsequence of the highest-degree homogeneous parts.
When these parts form a semi-regular sequence in nondecreasing degree order, the bound involves the $n$ smallest input degrees together with the largest one.
The other bound is expressed in terms of the saturation exponent with respect to the homogenizing variable.
Both are obtained by bounding the degree from which the Hilbert function of the quotient ring associated with the homogenized system is constant.

We then compare the Buchberger-like Gr\"obner basis computations for an affine system, its homogenization, and its highest-degree homogeneous parts.
The first degree fall is characterized by failure of injectivity of multiplication by the homogenizing variable.
Before that point, choices of S-pairs and reducers in any computation can be matched in the others, and reduction sequences, remainders, intermediate bases, and leading monomials correspond under specialization.
Cryptographic semi-regularity guarantees this correspondence until the step degree first reaches the degree of regularity.
At that degree, affine reduction steps that preserve the sugar degree lift to homogeneous ones, yielding upper bounds on the algorithmic solving degree for a computation starting directly from the affine input.
}

\section{Introduction}\label{sec:Intro}

Computing Gr\"obner bases is a fundamental task in computational commutative algebra and in solving polynomial systems.
Since Buchberger's original algorithm~\cite{Buchberger}, several major improvements, including the $F_4$~\cite{F4} and $F_5$~\cite{F5} algorithms, have been developed.
A central problem in both theory and practice is to estimate the complexity of these computations.
When the ideal generated by the input is zero-dimensional and a graded monomial ordering is used, particularly a degree reverse lexicographic (DRL) ordering, this complexity is commonly studied through various notions of \emph{solving degree}, which we review in Section~\ref{subsec:complexity}.
For homogeneous systems, degree growth can often be controlled by invariants of the generated ideal.
For affine systems, homogenization makes homogeneous methods available but leaves two distinct problems:
bounding degree growth in the homogenized computation and transferring homogeneous information to a computation starting directly from the affine input.
These are the two problems addressed in this paper.

Let $K$ be a field, and let $R=K[x_1,\ldots,x_n]$ be the polynomial ring in $n$ variables over $K$.
Throughout, we fix a DRL ordering $\prec$ on $R$.
For a polynomial $f\in R\smallsetminus \{0\}$, let $f^{\rm top}$ denote its homogeneous component of highest total degree, and let $f^h$ denote its homogenization in $R'=R[y]$ with respect to an extra variable $y$.
We denote by $\prec^h$ the homogenization of $\prec$ on $R'$, with $y$ as the smallest variable.
\textcolor{black}{Note that $\prec^h$ is also a DRL ordering.}
Let $\bm F=(f_1,\ldots,f_m)$ be a sequence of polynomials in $R\smallsetminus \{ 0 \}$, and put $F=\{f_1,\ldots,f_m\}$ and $d_i=\deg(f_i)$, where $\deg$ denotes total degree.
We write $\bm F^{\rm top}=(f_1^{\rm top},\ldots,f_m^{\rm top})$ and $\bm F^h=(f_1^h,\ldots,f_m^h)$, and use $F^{\rm top}$ and $F^h$ for the corresponding sets.
For a subset $H$ of $R$ or $R'$, we denote by $\langle H\rangle$ the ideal generated by $H$.
For a sequence $\bm H=(h_1,\ldots,h_r)$, we also write $\langle\bm H\rangle:=\langle h_1,\ldots,h_r\rangle$.
We write $d_{\rm reg}(I)$ for the degree of regularity of a homogeneous ideal $I$, and $\mathrm{max.GB.deg}_{\prec'}(I)$ for the maximum total degree of the elements in the reduced Gr\"obner basis of an ideal $I$ with respect to a monomial ordering $\prec'$; precise definitions are recalled in Section~\ref{sec:pre}.
When $I=\langle H\rangle$, we also write $\mathrm{max.GB.deg}_{\prec'}(H)$ for $\mathrm{max.GB.deg}_{\prec'}(I)$.

A general algebraic framework for studying solving degrees was developed by Caminata, Gorla, and their collaborators, who established rigorous relations among several notions of solving degree, the maximum Gr\"obner basis degree, and invariants from commutative algebra; see \cite{CG20,GMP22,CG23}.
A complementary line of research uses semi-regularity in the sense of Pardue~\cite{Pardue} and the weaker notion of cryptographic semi-regularity as tractable models for overdetermined polynomial systems.
The foundational work of Bardet, Faug\`ere, Salvy, and Yang described the Hilbert series and obtained asymptotic estimates for the degree of regularity under cryptographic semi-regularity; see \cite{Bar,BFS,BFSY}.
These estimates have played an important role in complexity analyses in multivariate cryptography.
Further work clarified the hypotheses under which such estimates give rigorous bounds for particular notions of solving degree and particular Gr\"obner basis algorithms; see \cite{BNDGMT21,GG23-2}.

The degree of regularity of $F^{\rm top}$ alone does not determine the degree growth in the computation for $F^h$, while relations among the final Gr\"obner bases do not determine the intermediate S-pairs and reductions in a computation starting from $F$.
In our previous work~\cite{KY}, we studied the Hilbert--Poincar\'e series of $R'/\langle F^h\rangle$ and compared the Gr\"obner bases of $\langle F\rangle$, $\langle F^h\rangle$, and $\langle F^{\rm top}\rangle$.
Building on those results, we address the two issues above by deriving new upper bounds for $\mathrm{max.GB.deg}_{\prec^h}(F^h)$ and by comparing, step by step, the Buchberger-like computations for $F$, $F^h$, and $F^{\rm top}$.

\paragraph{Degree bounds}
Our first group of results concerns degree growth in the homogenized system.
For homogeneous input, $\mathrm{max.GB.deg}_{\prec^h}(F^h)$ coincides with the two notions of solving degree reviewed in Section~\ref{subsec:complexity} that do not depend on a particular execution, and the relations in Section~\ref{subsec:knownbound} transfer a bound for this quantity to the corresponding solving degrees of $F$.
Lazard's classical bound provides a natural benchmark:
In the overdetermined case, its formula involves the $n+1$ largest input degrees.
Our first bound instead exploits a regular initial subsequence of $\bm F^{\rm top}$, while the second reflects saturation with respect to the homogenizing variable $y$.
The results are summarized as follows.

\begin{theor}[Theorem~\ref{thm:maxGB} and
Proposition~\ref{pro:dsat}]\label{thm:mainB}
Let $\bm F=(f_1,\ldots,f_m)$ be a sequence of nonconstant polynomials
in $R$, and put $d_i=\deg (f_i)$.
\begin{enumerate}
\item[(1)] Assume that $m>n$, that $1\le d_{n+1}\le\cdots\le d_m$, and that $(f_1^{\rm top},\ldots,f_n^{\rm top})$ is a regular sequence.
Then
\[
\mathrm{max.GB.deg}_{\prec^h}(F^h) \le d_1+\cdots+d_n+d_m-n.
\]
If, in addition,
\[
\max\left\{
d_m,\,
d_{\rm reg}
\bigl(\langle f_1^{\rm top},\ldots,f_{n+1}^{\rm top}\rangle\bigr)
\right\}
\le
\left\lfloor
\frac{d_1+\cdots+d_{n+1}-n-1}{2}
\right\rfloor+1,
\]
then
\[
\mathrm{max.GB.deg}_{\prec^h}(F^h)
\le
d_1+\cdots+d_n+d_{n+1}-n.
\]

\item[(2)] Put $D=d_{\rm reg}(\langle F^{\rm top}\rangle)$, and let $s$ be the least positive integer such that $(\langle F^h\rangle:y^s)=(\langle F^h\rangle:y^\infty)$, where $(I:y^\infty):=\bigcup_{k\geq1}(I:y^k)$ denotes the saturation of $I$ with respect to $y$.
Then $\mathrm{max.GB.deg}_{\prec^h}(F^h)\le D+s-1$.
\end{enumerate}
\end{theor}

The hypothesis in Theorem~\ref{thm:mainB}~{\rm (1)} on regularity is generic:
for fixed degrees over an infinite field, a generic choice of $f_1^{\rm top},\ldots,f_n^{\rm top}$ forms a regular sequence.
If $\bm F^{\rm top}$ is semi-regular, then its first $n$ elements form a regular sequence, and Theorem~\ref{thm:mainB}~{\rm (1)} applies after permuting only $f_{n+1},\ldots,f_m$ so that $d_{n+1}\leq\cdots\leq d_m$.
In particular, if the generators are indexed in nondecreasing degree order and $\bm F^{\rm top}$ is semi-regular in this order, the first bound involves the $n$ smallest input degrees together with the largest one, whereas Lazard's bound involves the $n+1$ largest input degrees.
Under the additional numerical condition in the theorem, $d_m$ is replaced by $d_{n+1}$, the $(n+1)$-st smallest input degree in this setting.
The second bound is system-dependent in a different sense.
Since $\langle F\rangle^h=(\langle F^h\rangle:y^\infty)$, the exponent $s$ measures the $y$-torsion separating the ideal generated by the homogenized input from the homogenization of the affine ideal.
If $\langle F^h\rangle$ is already saturated with respect to $y$, then $s=1$ and the bound reduces to $\mathrm{max.GB.deg}_{\prec^h}(F^h)\leq D$.
Thus this estimate complements the degree-based bounds above and connects the first group of results with the $y$-torsion governing degree falls below.

Both bounds are obtained by first estimating $\widetilde d_{\rm reg}(\langle F^h\rangle)$, defined in Section~\ref{subsec:ext} as the least degree from which the Hilbert function of $R'/\langle F^h\rangle$ is constant.
For Theorem~\ref{thm:mainB}~{\rm (1)}, this estimate follows from the vanishing of the first Koszul homology.
For Theorem~\ref{thm:mainB}~{\rm (2)}, in the nontrivial case $D<\infty$, it follows from an exact sequence associated with saturation by powers of $y$.
Together with a bound for $d_{\rm reg}(\langle F^{\rm top}\rangle)$, Hashemi's result converts these estimates into the stated bounds on the maximum Gr\"obner basis degree.

\paragraph{Degree falls and specialization}

Our second group of results concerns the Gr\"obner basis computation process itself.
The systems $F^{\rm top}$, $F^h$, and $F$ are related by the specializations $y=0$ and $y=1$.
The basic relation between an affine system and its highest-degree homogeneous parts before degree falls is familiar in the study of semi-regular systems and is implicit in several complexity analyses based on the degree of regularity, especially for matrix-$F_5$ type algorithms; see, e.g., \cite{Bar,BFS,BFSY,BFS-F5}.
Our contribution is to make this relation explicit at the level of individual steps of a Buchberger-like computation and determine how long the choices and the resulting polynomials can be made to correspond under specialization.
Here the step degree is the degree of the S-pair selected under the normal selection strategy; see Section~\ref{subsec:complexity}.

We also use this comparison to revisit earlier upper bounds on algorithmic solving degrees.
Tenti~\cite{Tenti} and Semaev--Tenti~\cite{ST} obtained such bounds over finite fields under assumptions involving the field equations.
In our previous work~\cite{KY}, analogous bounds over arbitrary fields were obtained by restarting the computation after dehomogenizing a partial Gr\"obner basis of the homogenized system.
The comparison developed here shows that this restart is unnecessary.
The main consequences are summarized as follows.

\begin{theor}\label{thm:mainC}
With notation as above, assume that $\bm F^{\rm top}$ is cryptographic semi-regular, and put $D=d_{\rm reg}(\langle F^{\rm top}\rangle)\in \mathbb{Z}_{>0}\cup \{\infty\}$.
\begin{enumerate}
\item[(1)]
If $D<\infty$, then, before the step degree first reaches $D$, corresponding S-pairs can be selected in the three Gr\"obner basis computations for $\langle F^{\rm top}\rangle$, $\langle F^h\rangle$, and $\langle F\rangle$.
After a reduction sequence is chosen in any one of them, corresponding sequences can be constructed in the other two;
reduction steps that specialize to the identity at $y=0$ are omitted.
The resulting remainders, intermediate bases, and leading monomials correspond under the specializations $y=0$ and $y=1$.
If $D=\infty$, the same statements hold through every finite step degree.
In particular, no degree fall occurs in the indicated part of the affine computation.

\item[(2)]
Assume in addition that $D<\infty$.
At the first stage at which the step degree reaches $D$, we carry out the affine computation using the sugar strategy and the reduction convention introduced in Definition~\ref{def:r(g1,g2)}.
For every selected affine S-pair of sugar degree $D$, the part of its reduction sequence that preserves the sugar degree lifts to a reduction sequence for the corresponding homogeneous S-pair.
Consequently, the leading monomials obtained by the end of sugar degree $D$ generate a monomial ideal containing every monomial of $R$ of degree $D$.

If $D\geq\max\{\deg(f):f\in F\}$, then the algorithmic solving degree
of the resulting Buchberger-like computation is at most $2D-1$.
If, after sugar degree $D$, the sugar strategy is stopped and the
polynomials obtained up to that degree are inter-reduced before the
computation is continued, then the solving degree in the strict sense
is at most $2D-2$.
\end{enumerate}
\end{theor}

Theorem~\ref{thm:mainC}~{\rm (1)} allows an S-pair or a reduction sequence to be chosen first in any one of the three computations; corresponding choices can then be made in the other two.
In particular, this correspondence is realized by the deterministic choices described in Remark~\ref{rem:deterministic}.
The mechanism behind the correspondence is the multiplication-by-$y$ map on $R'/\langle F^h\rangle$.
We define the $y$-injectivity threshold $d_{\rm inj}$ as the first positive degree in which this map is not injective and prove that it is exactly the first step degree at which a degree fall occurs in the affine computation.
If $d_{\rm inj}=\infty$, no degree fall occurs.
Cryptographic semi-regularity implies $D=d_{\rm reg}(\langle F^{\rm top}\rangle)\leq d_{\rm inj}$, but the converse does not hold; see Section~\ref{app:example-noncsr}.
Thus cryptographic semi-regularity is a sufficient condition for the correspondence below $D$, whereas injectivity of multiplication by $y$ is the underlying condition.
When the step degree first reaches $D$, a full correspondence of reduction sequences may fail.
Theorem~\ref{thm:mainC}~{\rm (2)} shows that the reduction steps preserving the sugar degree still lift from the affine computation to the homogeneous one.
This is enough to control all monomials of $R$ of degree $D$ and to apply the argument of~\cite[Lemma~4]{KY} directly to the computation starting from $F$.
The same construction can more generally be carried out at any finite threshold below which multiplication by $y$ is injective, in particular at $d_{\rm inj}$ when it is finite.

Beyond the numerical bounds, Theorem~\ref{thm:mainC} transfers degree-by-degree information from the homogeneous computations to a computation performed directly on the affine input.
Before the first degree fall, leading monomials, intermediate bases, and reduction outcomes may be analyzed on the homogeneous side without replacing $F$ by a dehomogenized partial Gr\"obner basis.
Because the governing condition is injectivity of multiplication by $y$, rather than semi-regularity itself, the same argument applies to systems outside the cryptographic semi-regular class whenever that injectivity can be established; see Section~\ref{app:example-noncsr}.
The theorem also permits arbitrary matched choices of S-pairs and reduction sequences within the stated setting, while Remark~\ref{rem:deterministic} gives one deterministic realization.

\section{Preliminaries and auxiliary results}\label{sec:pre}

In this section, we review several notions of solving degree and collect known degree bounds relevant to this paper.
In Section~\ref{subsec:ext}, we use the stabilization index of the Hilbert function, which agrees with the classical Hilbert regularity when finite, as a generalized degree of regularity.
We then relate this invariant to the maximum Gr\"obner basis degree and record a specialization property of Koszul homology that will also be used later.
We recall below the notions of semi-regularity, cryptographic semi-regularity, and their affine variants;
see also \cite[Section~2.1]{KY}, where concise descriptions and characterizations in terms of Hilbert series are given.

\paragraph{Notations and Terminology}
We continue to use the notation from Section~\ref{sec:Intro}.
In particular, let $R=K[x_1,\ldots,x_n]$ denote the polynomial ring in $n$ variables over a field $K$, let $F$ be a finite subset of $R\smallsetminus\{0\}$, and let $\prec$ be the fixed DRL ordering on $R$.
Unless otherwise stated, we follow \cite{Singular} for terminology and notation related to Gr\"obner bases.
For example, the total degree of a polynomial $f$ is denoted by $\deg(f)$. 
For a monomial ordering $\prec$, we denote by $\mathrm{LM}_{\prec}(f)$ (and $\mathrm{LT}_{\prec}(f)$) the leading monomial (and term) of a nonzero polynomial $f$ with respect to $\prec$. 
For a set $H$ of nonzero polynomials, we put $\operatorname{LM}_{\prec}(H):=\{\operatorname{LM}_{\prec}(h):h\in H\}$.
For an integer $d$, we also put $H_d:=\{h\in H:\deg(h)=d\}$ and $H_{<d}:=\{h\in H:\deg(h)<d\}$, and define $H_{\leq d}$ similarly.
When $d=\infty$, we use the convention $H_{<d}=H$.
For a set $H\subset R'=R[y]$ and $a\in K$, we write $H|_{y=a}:=\{h|_{y=a}:h\in H\}$; specialization of a sequence is understood componentwise.
For an ideal $I\subset R$, its homogenization in $R'$ is denoted by $I^h$.

For the definition of homogenization, see \cite[Chapter~2]{KR} or \cite[Section~3.1.1]{KY}.
For convenience, we recall the homogenization ordering $\prec^h$ on $R'=R[y]$ associated with $\prec$.
The ordering $\prec^h$ is a monomial ordering on $R'=R[y]$ defined as follows:
For two monomials $t_1=X_1y^a$ and $t_2=X_2y^b$, where $X_1$ and $X_2$ are monomials in $R$, we set $t_1\prec^h t_2$ if and only if {\rm (i)} $\deg(t_1)=\deg(X_1)+a<\deg(t_2)=\deg(X_2)+b$, or {\rm (ii)} $\deg(t_1)=\deg(t_2)$ and $X_1\prec X_2$.
\textcolor{black}{It is obvious that if $\prec$ is graded, then $\prec^h$ is also graded. 
Moreover, if $\prec$ is the DRL ordering with $x_n\prec\cdots \prec x_1$, then $\prec^h$ is also the DRL ordering with $y \prec x_n \prec \cdots \prec x_1$.}

For a graded $K$-algebra $A$ and an integer $d$, we denote by $A_d$ its homogeneous component of degree $d$.
For a non-negative integer $d$, we also put $R_{\leq d}:=\bigoplus_{j=0}^dR_j$.
For a finitely generated standard graded $K$-algebra $A=\bigoplus_{d\geq0}A_d$, we denote its Hilbert function, Hilbert series, and Hilbert polynomial by $\mathrm{HF}_A$, $\mathrm{HS}_A(z)$, and $\mathrm{HP}_A$, respectively, where $\mathrm{HF}_A(d)=\dim_K A_d$ and $\mathrm{HS}_A(z)=\sum_{d\geq0}\mathrm{HF}_A(d)z^d$.

For a commutative ring $A$, we denote by $\dim(A)$ its Krull dimension; see \cite[Section~10.2]{Eisen}.
For an ideal $I$ of a polynomial ring $S$, we say that $I$ is \emph{zero-dimensional} if $\dim(S/I)=0$; see also \cite[Chapter~IX]{CLO}.
Recall that, for a homogeneous ideal $I$ in a standard graded polynomial ring $S$, the degree of regularity $d_{\rm reg}(I)$ is defined as the smallest non-negative integer $d$ such that $I_d=S_d$, with the convention $d_{\rm reg}(I)=\infty$ if no such integer exists.
If $d_{\rm reg}(I)$ is finite, then it coincides with the Castelnuovo--Mumford regularity of $I$ as a graded module, denoted by $\mathrm{reg}(I)$; see \cite[Section~20.5]{Eisen} for a definition.

For a positive integer $d$, a homogeneous sequence $\bm H=(h_1,\ldots,h_m)\in R_{d_1}\times\cdots\times R_{d_m}$ is said to be \emph{$d$-regular} if $\mathrm{HS}_{R/\langle\bm H\rangle}(z) \equiv \frac{\prod_{j=1}^m(1-z^{d_j})}{(1-z)^n} \pmod{z^d}$.
We use this definition for every positive integer $d$, thereby extending its usual range $d\geq\max\{d_1,\ldots,d_m\}$.
The sequence $\bm H$ is said to be \emph{cryptographic semi-regular} if it is $d_{\rm reg}(\langle\bm H\rangle)$-regular;
when $d_{\rm reg}(\langle\bm H\rangle)=\infty$, this means that $\bm H$ is $d$-regular for every positive integer $d$.
For $1\leq i\leq m$, put $\bm H^{(i)}=(h_1,\ldots,h_i)$.
We say that $\bm H$ is \emph{semi-regular} if every initial subsequence $\bm H^{(i)}$ is cryptographic semi-regular.
By \cite[Propositions~1 and~3]{KY}, this is equivalent to semi-regularity in the sense of Pardue~\cite{Pardue}.
In particular, if $m\leq n$, then regularity, semi-regularity, and cryptographic semi-regularity are equivalent.
A sequence $\bm F$ of not necessarily homogeneous polynomials is said to be \emph{affine semi-regular} (resp.\ \emph{affine cryptographic semi-regular}) if $\bm F^{\rm top}$ is semi-regular (resp.\ cryptographic semi-regular).

For a graded module $M$ and an integer $a$, the notation $M(-a)$
denotes the internal degree shift defined by $M(-a)_d=M_{d-a}$.
If $C_\bullet$ is a complex of graded modules, then $C_\bullet(-a)$ denotes the complex obtained by applying this shift termwise: $(C_\bullet(-a))_i=C_i(-a)$.
Thus $(-a)$ changes only the internal grading and does not change the homological indices.

For a finite homogeneous sequence $\bm H=(h_1,\ldots,h_m)$ in a standard graded polynomial ring $S$, we denote by $K_\bullet(\bm H)$ its Koszul complex, regarded as a complex of graded free $S$-modules, with the grading convention of \cite[Definition~7.6.6]{Singular}.
Its boundary maps preserve the internal degree.
We denote its $i$-th homology by $H_i(K_\bullet(\bm H))$ and write $H_i(K_\bullet(\bm H))_d$ for its homogeneous component of internal degree $d$.
The ambient polynomial ring will be specified whenever it is not clear from the context.
For $d\geq\max\{d_1,\ldots,d_m\}$, Diem showed in \cite[Theorem~1]{Diem2} that $\bm H$ is $d$-regular if and only if $H_1(K_\bullet(\bm H))_t=0$ for every $t<d$.
The same equivalence remains valid for every positive integer $d$ under the convention above.
Indeed, deleting from $\bm H$ all elements of degree at least $d$ changes neither the Hilbert-series congruence modulo $z^d$ nor the Koszul complex in internal degrees below $d$.

\subsection{Solving degrees of Gr\"{o}bner basis computation}\label{subsec:complexity}

For Buchberger's algorithm and its variants with the normal selection strategy (see \cite[Section~II.10]{CLO}), S-pairs are selected according to the least common multiples of their leading monomials.
For an S-pair $\bm{s}=(g_1,g_2)$, we put $\delta_\prec(\bm{s}):=\deg(\operatorname{lcm}(\operatorname{LM}_\prec(g_1),\operatorname{LM}_\prec(g_2)))$ and call it the \emph{degree of the S-pair}.
The corresponding S-polynomial is
\begin{equation}\label{eq:S-poly}
S(g_1,g_2)
=
\frac{
\operatorname{lcm}
(\operatorname{LM}_\prec(g_1),\operatorname{LM}_\prec(g_2))
}{
\operatorname{LT}_\prec(g_1)
}
g_1
-
\frac{
\operatorname{lcm}
(\operatorname{LM}_\prec(g_1),\operatorname{LM}_\prec(g_2))
}{
\operatorname{LT}_\prec(g_2)
}
g_2.
\end{equation}
The two multiples on the right-hand side of \eqref{eq:S-poly} have total degree $\delta_\prec(\bm{s})$.

At each stage, the normal selection strategy processes an available S-pair whose least common multiple is minimal with respect to $\prec$.
Since $\prec$ is graded, the selected pair has minimum degree among all available S-pairs.
We call the degree of the selected S-pair the \emph{step degree}.
For homogeneous input, the sequence of step degrees is nondecreasing.
For inhomogeneous input, however, it may decrease after a reduction produces a polynomial of degree smaller than the degree of the selected S-pair.
The normal selection strategy is thus naturally adapted to computations with graded monomial orderings.

Degree growth is one of the main parameters used to estimate the complexity of Gr\"obner basis computation with respect to a graded monomial ordering.
The term \emph{solving degree} is used in the literature for several closely related, but not identical, degree parameters.
In this paper, we distinguish a solving degree attached to a particular execution from two \textcolor{black}{others} that do not depend on the choice of an execution.

Ding and Schmidt~\cite{Ding2013} call a computationally dominant step a \emph{solving step}, and call the highest degree of the polynomials involved in such a step the \emph{solving degree}.
For the theoretical degree bounds considered in this paper, we use the following closely related global formulation.

\begin{definition}[Solving degree I (algorithmic solving degree);
cf.\ {\cite[p.~36]{Ding2013}}]\label{def:solvedeg1}
Let $\mathcal A$ be a Gr\"obner basis algorithm applied to an input set $F$ with respect to a graded monomial ordering $\prec$.
The \emph{algorithmic solving degree} of $\mathcal A$ on $F$ is the highest total degree of the polynomials involved during the entire execution of $\mathcal A$.
We denote it by $\mathrm{sd}_{\prec}^{\mathcal A}(F)$.
\end{definition}

For the Buchberger-like algorithm considered in this paper, the algorithmic solving degree is the maximum of the highest step degree and the largest input degree.
Indeed, the input polynomials and the two multiples used to form each selected S-polynomial are regarded as polynomials involved in the computation.
The two multiples associated with an S-pair $\bm s$ have degree $\delta_\prec(\bm s)$, and reductions with respect to a graded monomial ordering do not increase total degree.
In particular, if the highest step degree is at least $\max\{\deg(f):f\in F\}$, then $\mathrm{sd}_{\prec}^{\mathcal A}(F)$ is equal to the highest step degree.

Note that an S-pair discarded by a criterion (cf.\ \cite[p.\ 116]{CLO}) before its S-polynomial is formed is not regarded as processed.
Accordingly, neither the degree of such a pair nor the corresponding multiples contributes to $\mathrm{sd}_{\prec}^{\mathcal A}(F)$.

\begin{remark}
The original notion of Ding and Schmidt is operational:
A solving step is selected according to running time or memory consumption.
It may therefore depend not only on the algorithm but also on its implementation and on the chosen cost measure.

Definition~\ref{def:solvedeg1}, by contrast, records the maximum degree reached during the whole computation.
It therefore gives an upper bound for the operational solving degree of Ding and Schmidt, and the two agree whenever a computationally dominant step occurs at the maximum degree reached.
\end{remark}

\begin{remark}[Solving degree in the strict sense]\label{rem:solvedeg1_str}
For a selected S-pair $\bm s=(g_1,g_2)$, the leading terms of the two multiples in \eqref{eq:S-poly} cancel, and hence one has $\deg (S(g_1,g_2))\leq\delta_\prec(\bm s)$.
Following the convention used in~\cite{Tenti,ST}, we call the highest
degree of the nonzero S-polynomials appearing during the computation the \emph{solving degree in the strict sense}.
Thus the strict solving degree records the actual degrees after the leading-term cancellations, rather than the degrees assigned to the selected S-pairs.
In particular, if $g_i^{\rm top}=\operatorname{LT}_\prec(g_i)$ for $i=1,2$ and $S(g_1,g_2)\neq0$, then $\deg (S(g_1,g_2))<\delta_\prec(\bm s)$.
\end{remark}

We next recall two notions of solving degree that do not depend on a particular execution.
The first, due to Caminata and Gorla~\cite{CG20}, is formulated in terms of Macaulay matrices.
For a fixed graded monomial ordering $\prec$ and a non-negative integer $d$, let $M_{\le d}(F)$ denote the degree-$d$ Macaulay matrix associated with $F$ and $\prec$ in the sense of \cite[Section~3]{CG20}.

\begin{definition}[Solving degree II~\cite{CG20}]\label{def:solvedeg2}
The solving degree of $F$ with respect to $\prec$ is defined as the smallest non-negative integer $d$ such that the reduced row echelon form of the degree-$d$ Macaulay matrix $M_{\le d}(F)$ yields a Gr\"obner basis of $F$ with respect to $\prec$.
We denote this degree by $\mathrm{sd}_{\prec}^{\mathrm{mac}}(F)$.
\end{definition}

The computation of the reduced row echelon form of $M_{\le d}(F)$ corresponds to the standard XL algorithm~\cite{XL}, which is based on an idea of Lazard~\cite{Lazard}.

The third definition is related to the last fall degree.
For each non-negative integer $d$, let $V_{F,d}$ be the smallest $K$-vector subspace of $R$ such that $F\cap R_{\leq d}\subset V_{F,d}$ and such that, whenever $f\in V_{F,d}$ and $g\in R$ satisfy $\deg(gf)\leq d$, one has $gf\in V_{F,d}$;
see \cite[Definition~1]{LFD2015}, \cite[Definition~2.1]{LFD2018}, and \cite[Definition~1.5]{CG23}.

\begin{definition}[Solving degree III]\label{def:solvedeg3}
The solving degree of $F$ with respect to the fixed graded monomial ordering $\prec$ is defined as the smallest non-negative integer $d$ such that the vector space $V_{F,d}$ contains a Gr\"obner basis of $F$ with respect to $\prec$.
We denote this degree by $\mathrm{sd}_{\prec}^{\mathrm{mut}}(F)$.
\end{definition}

The solving degree in Definition~\ref{def:solvedeg3} can be interpreted in terms of Macaulay matrices equipped with the so-called mutant strategy; this strategy computes a basis of $V_{F,d}$ for each $d$, and hence terminates precisely at the degree described in Definition~\ref{def:solvedeg3}; see \cite[Theorem~1]{GMP22}.
Algorithms employing this strategy, such as MutantXL~\cite{MutantXL} and MXL2~\cite{MXL2}, are often called \emph{mutant algorithms}; see also \cite{CG23,GG23-2,MXLsd}.
For this reason, we use the notation $\mathrm{sd}_{\prec}^{\mathrm{mut}}(F)$, where ``mut'' stands for the mutant strategy.

The two notions above correspond, respectively, to the solving degrees $\operatorname{solvdeg}^{s}$ and $\operatorname{solvdeg}^{m}$ for standard and mutant algorithms in~\cite{GG23-2}.
We use the superscripts ${\rm mac}$ and ${\rm mut}$ to indicate the definition using Macaulay matrices and the mutant strategy, respectively.

We next recall the notion of the maximum Gr\"obner basis degree, a purely algebraic invariant of the ideal that will be compared with the solving degrees.

\begin{definition}[Maximum Gr\"obner basis degree]\label{def:maxGBdeg}
Let $S$ be either $R$ or $R'$, let $I\subset S$ be an ideal, and let $\prec'$ be a monomial ordering on $S$.
The \emph{maximum Gr\"obner basis degree} of $I$ with respect to $\prec'$, denoted by $\mathrm{max.GB.deg}_{\prec'}(I)$, is defined as the maximum of the total degrees of the elements in the reduced Gr\"obner basis of $I$ with respect to $\prec'$.
If $I=\langle H\rangle$, we sometimes write $\mathrm{max.GB.deg}_{\prec'}(H)$ instead of $\mathrm{max.GB.deg}_{\prec'}(I)$.
\end{definition}

\begin{remark}\label{rem:inequality}
In a series of works (cf.\ \cite{CG20,BNDGMT21,GMP22,CG23,GG23-2}), Gorla et al.\ established precise relations between the solving degrees $\mathrm{sd}_{\prec}^{\mathrm{mac}}(F)$ and $\mathrm{sd}_{\prec}^{\mathrm{mut}}(F)$ and several algebraic invariants associated with $F$.
For any graded monomial ordering $\prec$, one has
\begin{equation}\label{eq:sdGB}
  \mathrm{max.GB.deg}_{\prec}(F)
  \le \mathrm{sd}_{\prec}^{\mathrm{mut}}(F)
  \le \mathrm{sd}_{\prec}^{\mathrm{mac}}(F).
\end{equation}
Moreover, if all elements of $F$ are homogeneous, then all three quantities coincide.
\end{remark}

\begin{remark}[On graded monomial orderings]
In practice, Gr\"obner bases with respect to graded orderings, especially DRL, can often be computed substantially faster than bases with respect to nongraded orderings such as the lexicographical ordering.
DRL is also central to the theoretical analysis used in this paper, since the relevant degree bounds, homogenization results, and comparisons with regularity are formulated for graded or DRL orderings.
For zero-dimensional ideals, a Gr\"obner basis can subsequently be converted to another ordering by the FGLM algorithm~\cite{FGLM}; more generally, one may use the Gr\"obner walk~\cite{GroebnerWalk}.
For these practical and theoretical reasons, we work primarily with graded orderings and, when needed, zero-dimensional ideals.
\end{remark}

\subsection{Known bounds on solving degrees and maximum Gr\"obner basis degrees}\label{subsec:knownbound}

Here we consider a polynomial set $F=\{f_1,\ldots,f_m\}\subset R\smallsetminus\{0\}$,
where $R=K[x_1,\ldots,x_n]$, and a graded monomial ordering $\prec$.

If $F$ consists of homogeneous polynomials, then all three quantities in \eqref{eq:sdGB} coincide, whence any upper bound on $\mathrm{max.GB.deg}_{\prec}(F)$ also bounds both of these solving degrees.
For example, $\mathrm{max.GB.deg}_{\prec}(F)$ can be bounded in terms of the Castelnuovo--Mumford regularity $\mathrm{reg}(\langle F \rangle)$ (cf.\ \cite[Corollary 2.5]{BS} and \cite[Proposition 4.15]{Hashemi2010}), under additional assumptions such as $\prec$ being a DRL ordering and $\langle F \rangle$ being in strong Noether position (see Remark \ref{rem:regGB} below).
The algorithmic solving degree $\mathrm{sd}_{\prec}^{\mathcal A}(F)$, however, depends on the particular execution, including the pair-discarding criteria and the stopping rule.
For the Buchberger-like computations used in our algorithmic results, the precise computational setting is specified at the beginning of Section~\ref{sec:app}.

\begin{remark}\label{rem:regGB}
Several results relate the maximum Gr\"obner basis degree to the Castelnuovo--Mumford regularity of a homogeneous ideal.
Throughout this remark, we assume that $\prec$ is the DRL ordering \textcolor{black}{with $x_n \prec \cdots \prec x_1$}.

Bayer and Stillman~\cite{BS} showed that, if a homogeneous ideal $I$ is given in generic coordinates, then $\mathrm{max.GB.deg}_{\prec}(I) \le \mathrm{reg}(I)$.
Hashemi later obtained related statements under structural assumptions on $I$.
In particular:
\begin{enumerate}
\item[(1)] By \cite[Theorem 4.17]{Hashemi2010}, if $I$ is in \emph{strong Noether position} (SNP), then
\begin{equation}\label{eq:Hashemi2010}
    \mathrm{max.GB.deg}_{\prec}(I) \le \mathrm{reg}(I)
= \overline{\mathrm{hilb}}(I),
\end{equation}
where $\overline{\mathrm{hilb}}(I) =
\max\left\{
\mathrm{hilb}(\mathrm{sec}(I,j))
\ \middle|\
0\leq j\leq\dim(R/I)
\right\}$ is the \emph{stabilized Hilbert regularity} introduced in
\cite{Hashemi2010}.
See \cite{Hashemi2010} for the definition of the sections $\mathrm{sec}(I,j)$ and Remark~\ref{rem:hilbreg} below for the definition of the Hilbert regularity $\mathrm{hilb}(-)$.

\item[(2)] By \cite[Theorem~4.1]{Hashemi2012}, if $I$ is in \emph{quasi stable position} (QSP) in the sense of \cite{Hashemi2012}, then
\begin{equation}\label{eq:Hashemi2012}
\mathrm{max.GB.deg}_{\prec}(I)=\mathrm{reg}(I).
\end{equation}
\end{enumerate}
We stress that the terminology is not uniform.
It follows from \cite[Lemma~3.1]{Hashemi2012} and \cite[Definition~3.1 and Proposition~3.2(v)]{HSS18} that every ideal in QSP in the sense of \cite{Hashemi2012} is in quasi-stable position in the sense of \cite{HSS18}.

The converse does not hold.
Indeed, for $2\leq i\leq n$ and $d\geq2$, the ideal $J_{i,d}:=\langle x_1^d,\ldots,x_i^d\rangle$ is quasi-stable in the sense of \cite{HSS18}, as is readily checked from \cite[Definition~3.1(i)]{HSS18}.
On the other hand, its minimal monomial generators form its reduced Gr\"obner basis, while $J_{i,d}$ is a complete intersection; hence $\mathrm{max.GB.deg}_{\prec}(J_{i,d})=d<i(d-1)+1=\mathrm{reg}(J_{i,d})$.
Therefore $J_{i,d}$ is not in QSP in the sense of \cite{Hashemi2012}, since QSP would imply equality by \eqref{eq:Hashemi2012}.
The special case $J_{2,2}\subset\mathbb Q[x_1,x_2]$ is explicitly noted in \cite[p.~1171]{Hashemi2012}.
Thus the QSP condition of \cite{Hashemi2012} is strictly stronger than quasi-stable position in the sense of \cite{HSS18}.

Finally, a homogeneous ideal is in quasi-stable position in the sense of \cite{HSS18} if and only if it is in strong Noether position in the sense of \cite{Hashemi2010};
see \cite[Definition~3.1 and Proposition~3.2(v)]{HSS18} and \cite[Definition~2.4, Theorem~2.19, and Corollary~2.20]{Hashemi2010}.
\end{remark}

In the remainder of this subsection, we consider the inhomogeneous case, i.e., $F$ contains at least one inhomogeneous polynomial.
In this case, the three quantities in \eqref{eq:sdGB} need not coincide, and it becomes difficult to estimate the solving degrees without \textcolor{black}{knowing} a Gr\"{o}bner basis of $\langle F \rangle$.

A straightforward way to bound the solving degrees in the inhomogeneous case is via homogenization.
For the rest of this subsection, we fix \textcolor{black}{a DRL ordering $\prec$ on $R$}, and denote by $\prec^h$ its homogenization on $R'=R[y]$.
Then, we have
\[
{\rm max.GB.deg}_{\prec}(F) \leq \mathrm{sd}^{\mathrm{mut}}_{\prec}(F) \leq \mathrm{sd}^{\mathrm{mac}}_{\prec}(F)
\leq  \mathrm{sd}^{\mathrm{mac}}_{\prec^h}(F^h)
= {\rm max.GB.deg}_{\prec^h}(F^h),
\]
where $F^h=\{f_1^h, \ldots, f_m^h\}$.
The inequality $\mathrm{sd}^{\mathrm{mac}}_{\prec}(F) \leq \mathrm{sd}^{\mathrm{mac}}_{\prec^h}(F^h)$ follows from the proof of \cite[Theorem~7]{CG20}, while the other inequalities and equality follow from Remark \ref{rem:inequality}.
We note that, without additional assumptions such as zero-dimensionality, the maximum Gr\"obner basis degree (and hence the solving degree) can grow doubly exponentially in the number of variables (see, e.g., Mayr--Ritscher~\cite{MR13}).
The bounds discussed below are therefore meaningful only under additional structural hypotheses.

For a zero-dimensional affine ideal satisfying some hypotheses below, Lazard's classical bound applies to $\mathrm{max.GB.deg}_{\prec^h}(F^h)$.

\begin{theorem}[{Lazard's bound}; {\cite[Theorem 2]{Lazard}}, {\cite[Th\'{e}or\`{e}m 3.3]{Lazard81}}]\label{thm:Lazard_zero_dim}
With notation as above, we assume that the number of projective zeros of $F^h$ is finite (and therefore $m \geq n$), and that $f_1^h=\cdots = f_m^h=0$ has no non-trivial solution over the algebraic closure $\overline{K}$ with $y=0$, i.e., $F^{\rm top}$ has no solution in $\overline{K}^n$ other than $(0,\ldots , 0)$.
Then, supposing also that $d_1 \geq \cdots \geq d_m$, where $d_i:=\deg(f_i)$ for $1\leq i\leq m$, we have
\begin{equation}\label{eq:Lazard0} \mathrm{max.GB.deg}_{\prec^h}(F^h) \leq d_1 + \cdots + d_{\ell} - \ell + 1
\end{equation}
with $\ell := \min\{ m,n+1 \}$.
\end{theorem}

For a detailed exposition of Lazard's argument and related degree bounds, see also \cite{C14}.
One of the key ingredients in the proof of Theorem~\ref{thm:Lazard_zero_dim} is the following lemma.

\begin{lemma}[cf.\ \cite{Lazard}]\label{lem:Lazard}
Let $I$ be a homogeneous ideal of $R'=R[y]$, and let $d_0$ be a positive integer satisfying the following two properties:
\begin{enumerate}
    \item[(1)] The multiplication-by-$y$ map $(R'/I)_{d_0-1} \longrightarrow (R'/I)_{d_0}$ is surjective.
    \item[(2)] For any integer $d \ge d_0$, the multiplication-by-$y$ map $(R'/I)_{d} \longrightarrow (R'/I)_{d+1}$ is injective.
\end{enumerate}
Then $\mathrm{max.GB.deg}_{\prec'}(I) \le d_0$, where $\prec'$ is the homogenization (with respect to $y$) of any graded monomial ordering on $R$.
\end{lemma}

Lazard showed that one can take $d_0=d_1+\cdots+d_{\ell}-\ell+1$ in Lemma~\ref{lem:Lazard} to obtain Theorem~\ref{thm:Lazard_zero_dim}, where we take $I=\langle F^h\rangle$ and $\prec'=\prec^h$ respectively.
The bound in \eqref{eq:Lazard0}, often referred to as the \emph{Macaulay bound}, yields an upper bound for $\mathrm{sd}_{\prec}^{\rm mac}(F)$ via homogenization, and hence also for $\mathrm{sd}_{\prec}^{\rm mut}(F)$.
Whether the same bound applies to an algorithmic solving degree or not depends on the computational setting.
The Buchberger-like setting used for our algorithmic results is specified in Section~\ref{sec:app}.

\paragraph{Bounds via the top homogeneous system}
For the maximum Gr\"obner basis degree of $\langle F\rangle$, if
$R/\langle F^{\rm top}\rangle$ is Artinian, 
{where $F^{\rm top}=\{f_1^{\rm top},\ldots,f_m^{\rm top}\}$}, then one has
\begin{equation}\label{eq:crypto}
\mathrm{max.GB.deg}_{\prec'}(F) \leq d_{\rm reg}(\langle F^{\rm top} \rangle)
\end{equation}
for any graded monomial ordering $\prec'$ on $R$; see~\cite[Remark~15]{CG20} or \cite[Lemma 4]{KY} for a proof.
When $d_{\rm reg}(\langle F^{\rm top}\rangle)<\infty$, it coincides with the Castelnuovo--Mumford regularity $\mathrm{reg}(\langle F^{\rm top}\rangle)$ of the homogeneous ideal $\langle F^{\rm top}\rangle$.

\begin{remark}
Note that both $d_{\rm reg}(\langle F^{\rm top}\rangle)$ and $\mathrm{sd}_{\prec}^{\rm mac}(F)$ are always greater than or equal to $\mathrm{max.GB.deg}_{\prec}(F)$.
In the cryptographic literature, the degree of regularity $d_{\rm reg}(\langle F^{\rm top}\rangle)$ (or related quantities such as the first fall degree) is often used as a proxy or heuristic upper bound for the solving degree, motivated by the inequality~\eqref{eq:crypto} and by its good predictive behavior in many structured instances.
On the other hand, it is pointed out in~\cite{BNDGMT21,CG20,CG23} by explicit examples that the degree of regularity does \emph{not} control the solving degrees in general:
Even when $\bm{F}=(f_1,\ldots,f_m)$ is an affine (cryptographic) semi-regular sequence (that is, $\bm{F}^{\rm top}$ is (cryptographic) semi-regular), the quantities $\mathrm{sd}_{\prec}^{\rm mac}(F)$ and $\mathrm{sd}_{\prec}^{\rm mut}(F)$ can be strictly larger than $d_{\rm reg}(\langle F^{\rm top}\rangle)$.
\end{remark}

As for $\mathrm{sd}_{\prec'}^{\rm mut}(F)$, Caminata and Gorla proved in~\cite{CG23} that it coincides with the \emph{last fall degree} $d_F$, introduced by Huang, Kosters, and Yeo~\cite[Definition 2]{LFD2015} (see also \cite[Definition 2.5]{LFD2018}) whenever $d_F > {\rm max.GB.deg}_{\prec'}(F)$.
In particular, if $d_{\rm reg}(\langle F^{\rm top}\rangle) < d_F$ and $d_F > {\rm max.GB.deg}_{\prec'}(F)$, then $d_{\rm reg}(\langle F^{\rm top}\rangle)<\mathrm{sd}_{\prec'}^{\rm mut}(F)$.
Hence $d_{\rm reg}(\langle F^{\rm top}\rangle)$ cannot provide an upper bound for $\mathrm{sd}_{\prec'}^{\rm mut}(F)$.

Nevertheless, a refined bound is available for the mutant strategy.
Recently, Salizzoni~\cite[Theorem~1.1]{MXLsd} proved that
\[
\mathrm{sd}_{\prec'}^{\rm mut}(F) \le d_{\rm reg}(\langle F^{\rm top}\rangle)+1,
\]
under the assumption that $d_{\rm reg}(\langle F^{\rm top}\rangle)\ge \max\{\deg(f): f\in F\}$.
Thus, although $d_{\rm reg}(\langle F^{\rm top}\rangle)$ alone does not bound the solving degrees in general, it still yields an upper bound for $\mathrm{sd}_{\prec'}^{\rm mut}(F)$ up to an additive constant of one.

\paragraph{Bounds for Buchberger-like algorithms with special settings}
In the case where $D=d_{\rm reg}(\langle F^{\rm top}\rangle)$ is finite, one can derive upper bounds on $\mathrm{sd}_{\prec}^{\mathcal{A}}(F)$ for certain Buchberger-like algorithms based on S-polynomial computation (including variants based on sub-Macaulay matrix computation).
These algorithms proceed in two phases:
The computation is first carried out up to a fixed degree bound and then continued using the intermediate polynomials obtained so far.
Semaev and Tenti \cite{ST} considered the case where $K=\mathbb{F}_q$ and $F$ contains the field equations $x_i^q-x_i$ for all $1\le i\le n$. 
In the first phase of their algorithm, a generating set $B$ of $\langle F\rangle$ is constructed by row-echelon form computation of sub-Macaulay matrices so that every monomial of degree at least $D$ is divisible by the leading monomial ${\rm LM}_{\prec}(g)$ for some $g\in B$.
The second phase then applies a Buchberger-like algorithm starting from $B$ as its input set.

\begin{theorem}[{\cite[Theorem 2.1]{ST}, \cite[Theorem 3.65 \& Corollary 3.67]{Tenti}}]\label{thm:tenti}
    With notation as above, assume that $K=\mathbb{F}_q$, and that $F$ contains $x_i^q - x_i$ for all $1 \leq i \leq n$.
    If $ D :=d_{\rm reg}(\langle F^{\rm top} \rangle) \geq \max \{ \deg (f) : f \in F \}$, then there exists a Buchberger-like algorithm $\mathcal{A}$ (with two phases) to compute the reduced Gr\"{o}bner basis of $F$ (based on S-polynomials) such that $\mathrm{sd}^{\mathcal{A}}_{\prec}(F) \leq 2 D  - 1$. 
    Furthermore, in the strict sense of Remark~\ref{rem:solvedeg1_str}, $\mathrm{sd}^{\mathcal{A}}_{\prec}(F) \leq 2 D  - 2$. 
\end{theorem}

Independently of Semaev and Tenti's work, analogous bounds over arbitrary fields and without assuming field equations were also obtained in our previous work~\cite[Lemma~4]{KY}.
The argument in~\cite{KY} was formulated under the additional assumption that $\bm F^{\rm top}$ is cryptographic semi-regular and considered another Buchberger-like algorithm with two phases.
In Section~\ref{subsec:overD} below, we revisit this result and show that similar bounds remain valid for Buchberger-like algorithms using the sugar strategy together with the reduction convention introduced in Definition~\ref{def:r(g1,g2)} below; see Lemma~\ref{lem:sd}.

\subsection{Hilbert regularity as a generalized degree of regularity}\label{subsec:ext}

Recall that $d_{\rm reg}(I)=\infty$ whenever the quotient by a homogeneous ideal $I$ is not Artinian.
In particular, this occurs for $\langle F^h\rangle$ when $F^h$ has a nontrivial projective zero over the algebraic closure $\overline K$.
If $F^h$ has only finitely many projective zeros---equivalently, if $R'/\langle F^h\rangle$ has Krull dimension at most one---then the Hilbert function of $R'/\langle F^h\rangle$ is eventually constant.
We use its stabilization index as a replacement for $d_{\rm reg}$ in the degree bounds below.

\begin{definition}[Generalized degree of regularity]\label{def:new}
Let $S$ be a standard graded polynomial ring over $K$, and let $I\subset S$ be a homogeneous ideal.
We define the \emph{generalized degree of regularity} of $I$ by
\[
\widetilde d_{\rm reg}(I)
:= \inf\left\{ d_0\in\mathbb Z_{\geq0} \ \middle|\ \dim_K(S/I)_d=\dim_K(S/I)_{d_0}
\text{ for all }d\geq d_0\right\},
\]
with the convention that $\inf\emptyset=\infty$.
\end{definition}

\begin{remark}\label{rem:hilbreg}
The \emph{index of regularity}, also called the \emph{Hilbert regularity} and often denoted by $\mathrm{hilb}(I)$, is the smallest integer $i\geq0$ such that ${\rm HF}_{S/I}(d)={\rm HP}_{S/I}(d)$ for all $d\geq i$.

Since the Hilbert function agrees with the Hilbert polynomial in all sufficiently large degrees, $\widetilde d_{\rm reg}(I)<\infty$ if and only if ${\rm HP}_{S/I}$ is constant, and in this case $\widetilde d_{\rm reg}(I)=\mathrm{hilb}(I)$.
In particular, if $S/I$ is Artinian, then ${\rm HP}_{S/I}=0$ and $\widetilde d_{\rm reg}(I)=d_{\rm reg}(I)=\mathrm{hilb}(I)$.
If $S/I$ is not Artinian, then $d_{\rm reg}(I)=\infty$, whereas $\widetilde d_{\rm reg}(I)$ may still be finite.
For instance, for $I=\langle y\rangle\subset K[x,y]$, one has $d_{\rm reg}(I)=\infty$ but $\widetilde d_{\rm reg}(I)=0$.

This identification allows us to apply bounds involving Hilbert regularity due to Hashemi and collaborators to the maximum Gr\"obner basis degree; see Proposition~\ref{prop:new} below.
\end{remark}

We record the following standard specialization property of Koszul homology.
It makes explicit the homological identification underlying the tensor-product argument in the proof of \cite[Theorem~7]{KY}.

\begin{lemma}\label{lem:koszul-specialization}
Let $\bm H=(h_1,\ldots,h_m)$ be a sequence of homogeneous polynomials in $R'=R[y]$, and put $\bm{H}'=(h_1,\ldots,h_m,y)$.
Then, for every $i\geq0$, specialization at $y=0$ induces a degree-preserving isomorphism of graded $K$-vector spaces $H_i\bigl(K_\bullet(\bm H')\bigr)\cong H_i\bigl(K_\bullet(\bm H|_{y=0})\bigr)$, where the Koszul complex on the left is taken over $R'$ and that on the right is taken over $R$.
\end{lemma}

\begin{proof}
Put $C_\bullet:=K_\bullet(\bm H)$ and $C'_\bullet:=K_\bullet(\bm H')$, where both Koszul complexes are taken over $R'$, and let $\partial_i:C_i\longrightarrow C_{i-1}$ and $\partial_i':C_i'\longrightarrow C_{i-1}'$ denote their boundary maps.
Regarding $C'_\bullet$ as the mapping cone of $C_\bullet(-1)\xrightarrow{\times y}C_\bullet$, we write $C_i'=C_i\oplus C_{i-1}(-1)$ and choose the signs so that
\[
\partial_i'(u,v)
=
\bigl(\partial_i(u)+yv,-\partial_{i-1}(v)\bigr).
\]
For every $i$, put $\overline C_i:=C_i/yC_i$.
Since $\partial_i(yC_i)\subseteq yC_{i-1}$, the boundary maps of
$C_\bullet$ induce maps $\overline\partial_i: \overline C_i\longrightarrow\overline C_{i-1}$, and hence the graded modules $\overline C_i$ form a complex $\overline C_\bullet$.
Define $\Phi_i:C_i'\longrightarrow\overline C_i$ by $\Phi_i(u,v)=u\bmod yC_i$.
Then the maps $\Phi_i$ are surjective and form a morphism of complexes, because
\[
\begin{aligned}
\Phi_{i-1}\bigl(\partial_i'(u,v)\bigr)=
\partial_i(u)+yv\bmod yC_{i-1}=
\partial_i(u)\bmod yC_{i-1}=
\overline\partial_i\bigl(\Phi_i(u,v)\bigr).
\end{aligned}
\]

We next show that the kernel complex of $\Phi$, denoted by $\ker(\Phi_\bullet)$, has zero homology.
Since each $C_i$ is a free $R'$-module, multiplication by $y$ on $C_i$ is injective.
Thus every element of $\ker(\Phi_i)$ can be written uniquely as $(yu,v)$ with $u\in C_i(-1)$ and $v\in C_{i-1}(-1)$.
Let
\[
\psi_i: C_i(-1)\oplus C_{i-1}(-1) \longrightarrow \ker(\Phi_i)
\ ;\ (u,v)\longmapsto(yu,v).
\]
On $C_i(-1)\oplus C_{i-1}(-1)$, let
\[
\delta_i(u,v) = \bigl(\partial_i(u)+v,-\partial_{i-1}(v)\bigr),
\]
which is the boundary map of the mapping cone of the identity map on $C_\bullet(-1)$.
A direct computation gives $\partial_i'\circ\psi_i=\psi_{i-1}\circ\delta_i$.
Hence the maps $\psi_i$ define an isomorphism of complexes between the mapping cone of the identity map on $C_\bullet(-1)$ and $\ker(\Phi_\bullet)$.
If $\delta_i(u,v)=0$, then $\partial_i(u)+v=0$, and hence $v=-\partial_i(u)$.
Therefore, $(u,v)=\bigl(u,-\partial_i(u)\bigr)=\delta_{i+1}(0,u)$.
Thus every cycle is a boundary, and consequently $H_i(\ker(\Phi_\bullet))=0$ for every $i$.

Therefore, the long exact sequence in homology associated with the short exact sequence of complexes $0\longrightarrow \ker(\Phi_\bullet) \longrightarrow C_\bullet' \xrightarrow{\ \Phi\ } \overline C_\bullet \longrightarrow0$ gives degree-preserving isomorphisms $H_i(C_\bullet') \cong H_i(\overline C_\bullet)$ for every $i$.

Finally, reduction modulo $y$ replaces each $h_j$ by
$h_j|_{y=0}$, and therefore gives an isomorphism of graded complexes $\overline C_\bullet \cong K_\bullet(\bm H|_{y=0})$, where the latter Koszul complex is taken over $R$.
The assertion follows.
\end{proof}

The following proposition provides an upper bound on the maximum Gr\"obner basis degree, and hence on the Macaulay and mutant solving degrees, of a homogeneous system $H$.
It will later be applied to $F^h$.

\begin{proposition}\label{prop:new}
Let $H = \{ h_1, \ldots, h_m \}$ be a set of homogeneous polynomials in $R' \smallsetminus K$ with $R'=R[y]$, and put $\bm{H}:=(h_1,\ldots,h_m) \in (R' \smallsetminus K)^m$.
Assume that $R/\langle H|_{y=0}\rangle$ is Artinian.
Then we have the following:
\begin{enumerate}
    \item[(1)] Assume that no element of $H$ is divisible by $y$, and that $\bm H|_{y=0}$ is cryptographic semi-regular.
    Then $\widetilde{d}_{\rm reg}( \langle H \rangle_{R'}) \ge d_{\rm reg}(\langle H|_{y=0} \rangle_R)-1$.

    \item[(2)] For any DRL ordering $\prec'$ on $R'$ such that $y$ is the least among the variables, one has $\mathrm{max.GB.deg}_{\prec'} (H ) \leq \max \{ {d}_{\rm reg}(\langle H|_{y=0} \rangle_{R}), \widetilde{d}_{\rm reg}( \langle H \rangle_{R'}) \}$. 
    This inequality is an equality if the monomial ideal $\langle {\rm LM}_{\prec'}(\langle H\rangle)\rangle $ is weakly reverse lexicographic;
    see \cite[Section~4]{Pardue} for the definition of a weakly reverse lexicographic monomial ideal.
\end{enumerate}
\end{proposition}

\begin{proof}
Put $D:=d_{\rm reg}(\langle H|_{y=0} \rangle_R)$ and $D':=\widetilde{d}_{\rm reg}(\langle H \rangle_{R'})$.
Let $K_\bullet$ and $K'_\bullet$ denote the Koszul complexes on $(h_1,\ldots,h_m)$ and $(h_1,\ldots,h_m,y)$ over $R'$, respectively.
Regarding $K'_\bullet$ as the mapping cone of $K_\bullet(-1)\xrightarrow{\times y}K_\bullet$, we obtain the exact sequence
\begin{equation}\label{eq:exact}
{
	\xymatrix{
	  H_{1}(K_{\bullet}')_d \ar[r] & H_{0}(K_{\bullet})_{d-1} \ar[r]^{\times y} & H_{0}(K_{\bullet})_d \ar[r] & H_{0}(K_{\bullet}')_d \ar[r] & 0    \\
	}}
\end{equation}
for each $d$.

We first prove (1).
Assume that no element of $H$ is divisible by $y$, and that
$\bm H|_{y=0}$ is cryptographic semi-regular.
By Lemma~\ref{lem:koszul-specialization}, for every $i$ and every $d$ one has $H_i(K_\bullet')_d \cong H_i\bigl(K_\bullet(\bm H|_{y=0})\bigr)_d$.
Hence, for every $d<D$, it follows from the definition of $d_{\rm reg}$ that $H_0(K_\bullet')_d\cong (R/\langle H|_{y=0}\rangle)_d \neq0$.
Moreover, the cryptographic semi-regularity of $\bm H|_{y=0}$ and \cite[Theorem~1]{Diem2} give $H_1(K_\bullet')_d=0$ for every $d<D$.
Hence exactness of~\eqref{eq:exact} implies $\dim_K(R'/\langle H\rangle)_{d-1}< \dim_K(R'/\langle H\rangle)_d$ for every $d<D$.
Therefore, we have $D' \geq D-1$, which proves (1).

We next prove (2).
Again by the definition of $d_{\rm reg}$, we have $H_0(K_{\bullet}')_d =0$ for all $d \geq D$.
Thus, for such $d$, exactness of~\eqref{eq:exact} shows that multiplication by $y$ induces a surjection $H_0(K_\bullet)_{d-1}\to H_0(K_\bullet)_d$, and therefore $\dim_K(R'/\langle H \rangle)_{d-1} \geq \dim_K (R'/\langle H \rangle)_{d}$ holds. 
Thus, the non-negative integral sequence $\left(\dim_K(R'/\langle H\rangle)_d\right)_{d\geq D}$ is weakly decreasing, from which we have $D'<\infty$, equivalently the Krull dimension of $R'/\langle H \rangle$ is at most one.
Moreover, the homogeneous ideal $\langle H\rangle$ is in Noether position.
Indeed, if $\dim(R'/\langle H\rangle)=1$, then $R'/(\langle H\rangle+\langle y\rangle)$ is Artinian, so that $R'/\langle H\rangle$ is finite over $K[y]$.
If $\dim(R'/\langle H\rangle)=0$, Noether position is automatic.

Since $\dim(R'/\langle H\rangle)\leq1$, Noether position is equivalent to quasi-stable position in the sense of \cite{HSS18}; see \cite[the paragraph following Proposition~5.7]{HS}.
By Remark~\ref{rem:regGB}, quasi-stable position in this sense is equivalent to strong Noether position.
Hence $\langle H\rangle$ is in strong Noether position, and
\eqref{eq:Hashemi2010} gives $\mathrm{max.GB.deg}_{\prec'}(H) \leq\mathrm{reg}(\langle H\rangle)=\overline{\mathrm{hilb}}(\langle H\rangle)$.
If $\dim(R'/\langle H\rangle)=1$, then $\mathrm{sec}(\langle H\rangle,0)=\langle H\rangle$ and $\mathrm{sec}(\langle H\rangle,1)=\langle H\rangle+\langle y\rangle$, so that $\overline{\mathrm{hilb}}(\langle H \rangle)=\max\{D',D\}$ by Remark~\ref{rem:hilbreg}.
If $\dim(R'/\langle H\rangle)=0$, then $\overline{\mathrm{hilb}}(\langle H \rangle)=D'$.
Moreover, $R'/(\langle H\rangle+\langle y\rangle)$ is a graded quotient of $R'/\langle H\rangle$ and is isomorphic to $R/\langle H|_{y=0}\rangle$; hence $D\leq D'$.
Thus $\overline{\mathrm{hilb}}(\langle H \rangle)=\max\{D,D'\}$ in either case, proving {\rm (2)}.

The equality under the weakly reverse lexicographic assumption is proved by applying general structural facts as follows:
Assume that the initial ideal $J:=\langle \mathrm{LM}_{\prec'}(\langle H\rangle)\rangle$ is weakly reverse lexicographic.
Then it follows directly from the definition of weakly reverse lexicographic ideals that $J$ is strongly stable; see also \cite[Definition 3.1 (iii)]{HSS18} for the definition of a strongly stable ideal.
In particular, $J$ satisfies the QSP criterion in \cite[Corollary~3.1]{Hashemi2012}, and hence $\langle H \rangle$ is in QSP in the sense of \cite{Hashemi2012}.
Therefore~\eqref{eq:Hashemi2012} gives $\mathrm{max.GB.deg}_{\prec'}(H)=\mathrm{reg}(\langle H \rangle)$.
Since the preceding argument shows that $\mathrm{reg}(\langle H \rangle)=\max\{D,D'\}$, equality holds in {\rm (2)}.
\end{proof}

\begin{remark}
Although Proposition~\ref{prop:new} is stated for homogeneous systems $H\subset R[y]$, it applies in particular to any affine polynomial system $F\subset R$ via homogenization.
Indeed, taking $H=F^h$, one has $H|_{y=0}=F^{\rm top}$.
Hence the assumptions of the proposition reduce to conditions on the top homogeneous system $F^{\rm top}$.
The resulting bound controls $\mathrm{max.GB.deg}_{\prec^h}(F^h)$ and hence, via homogenization, $\mathrm{sd}_{\prec}^{\rm mac}(F)$.
By~\eqref{eq:sdGB}, it also bounds $\mathrm{sd}_{\prec}^{\rm mut}(F)$.
\end{remark}

\section{New bounds on maximum Gr\"{o}bner basis degree}\label{sec:3}

In this section, we prove new upper bounds on the maximum Gr\"obner basis degree, using the notation introduced in Sections~\ref{sec:Intro} and~\ref{sec:pre}.
Throughout, $R=K[x_1,\ldots,x_n]$ is equipped with the fixed DRL ordering $\prec$, and $\prec^h$ denotes its homogenization on $R'=R[y]$.
For a sequence $\bm F=(f_1,\ldots,f_m)$ of non-constant polynomials in $R$, we put $d_j=\deg(f_j)$, $\bm F^{(k)}=(f_1,\ldots,f_k)$, and $F^{(k)}=\{f_1,\ldots,f_k\}$ for $k\le m$, and retain the notation $F$, $F^{\rm top}$, and $\bm F^h$ from Section \ref{sec:Intro}.

\subsection{New bounds under a regularity assumption}\label{sec:main}

Throughout the rest of this paper, for positive integers $n$ and $d_1,\ldots,d_m$ with $m \geq n$, we set
\begin{equation}\label{eq:newbound0}
  D^{(n;d_1,\ldots,d_m)}
  := \deg\!\left(\left[\frac{\prod_{j=1}^m (1-z^{d_j})}{(1-z)^n}\right]\right) + 1 ,
\end{equation}
where $[\cdot]$ means truncating a formal power series in $z$ over $\mathbb{Z}$ after the last consecutive positive coefficient.
Note that $D^{(n;d_1,\ldots,d_m)}$ does not depend on the order of $d_1,\ldots,d_m$, and $D^{(n+1;d_1,\ldots,d_m,1)} = D^{(n;d_1,\ldots,d_m)}$ for any $(n;d_1,\ldots,d_m)$.
It is known (see \cite[Lemma 2.5]{MM}) that
\[
D^{(n;d_1,\ldots,d_n)}=d_1+\cdots +d_n-n+1, \ 
D^{(n;d_1,\ldots,d_{n+1})}=\left\lfloor \frac{d_1+\cdots+d_{n+1}-n-1}{2}\right\rfloor+1.
\]

\begin{lemma}\label{lem:key}
Let $h_1,\ldots,h_m$ be non-constant homogeneous polynomials in $R$ of degrees $d_1,\ldots,d_m$.
Put $\bm{H} = (h_1,\ldots,h_m)$ and $H = \{ h_1,\ldots,h_m\}$.
Assume that $m>n$ and $1\le d_{n+1}\le\cdots\le d_m$ and that $\bm{H}^{(n)}=(h_1,\ldots,h_n)$ forms a regular sequence in $R$.
Then $H_1(K_\bullet(\bm{H}))_d=0$ for all $d\ge D^{(n;d_1,\ldots,d_n)}+d_m$.
Moreover, if $\max\{d_m, d_{\rm reg}(\langle H^{(n+1)}\rangle)\}\le D^{(n;d_1,\ldots,d_{n+1})}$ with $H^{(n+1)}=\{h_1,\ldots,h_{n+1}\}$, then $H_1(K_\bullet(\bm{H}))_d=0$ for all $d\ge D^{(n;d_1,\ldots,d_{n})}+d_{n+1}$.
\end{lemma}

\begin{proof}
For each $2\le j\le m$, we regard $K_\bullet(\bm{H}^{(j)})$ as the mapping cone of
\[
K_\bullet(\bm{H}^{(j-1)})(-d_j)\xrightarrow{\ \times h_j\ } K_\bullet(\bm{H}^{(j-1)}),
\]
as in~\cite{Diem2}.
This yields a long exact sequence of graded $K$-vector spaces; in particular, taking $j=m$ we obtain
\begin{equation}\label{eq:H1exact}
\xymatrix{
H_{1}(K_\bullet(\bm{H}^{(m-1)}))_{d} \ar[r] &
H_{1}(K_\bullet(\bm{H}))_{d} \ar[r] &
H_{0}(K_\bullet(\bm{H}^{(m-1)}))_{d-d_{m}} .
}
\end{equation}
We prove the assertions by induction on $m$.

\medskip\noindent
\emph{Case $m=n+1$.}
Since $\bm{H}^{(n)}$ is regular, we have $H_1(K_{\bullet}(\bm{H}^{(n)}))_d=0$ for all $d$ and $d_{\mathrm{reg}}(\langle H^{(n)}\rangle)=D^{(n;d_1,\ldots,d_n)}$.
Thus, it follows that $H_0(K_{\bullet}(\bm{H}^{(n)}))_{\ell}=(R/\langle H^{(n)}\rangle)_{\ell}=0$ for all $\ell \ge D^{(n;d_1,\ldots,d_n)}$.
Considering the sequence \eqref{eq:H1exact} for $m=n+1$, we obtain $H_1(K_{\bullet}(\bm{H}^{(n+1)}))_d=0$ for all $d\ge D^{(n;d_1,\ldots,d_n)}+d_{n+1}$, which proves the first assertion for $m=n+1$ (and implies the second assertion in this case).

\medskip\noindent
\emph{Induction step for the first assertion.}
Assume now that $m>n+1$ and that the first assertion holds for $m-1$.
Let $d\ge D^{(n;d_1,\ldots,d_n)}+d_m$.
Since $d_m\ge d_{m-1}$, we have $d\ge D^{(n;d_1,\ldots,d_n)}+d_{m-1}$, and thus $H_1(K_\bullet(\bm{H}^{(m-1)}))_d=0$ by the induction hypothesis.
It remains to show $H_0(K_\bullet(\bm{H}^{(m-1)}))_{d-d_m}=0$.
Note that $\langle {H}^{(n)}\rangle\subset \langle{H}^{(m-1)}\rangle$, hence $d_{\rm reg}(\langle {H}^{(m-1)}\rangle)\le d_{\rm reg}(\langle {H}^{(n)}\rangle)=D^{(n;d_1,\ldots,d_n)}$.
Since $d-d_m\ge D^{(n;d_1,\ldots,d_n)}$, we have $H_0(K_{\bullet}(\bm{H}^{(m-1)}))_{d-d_m}=(R/\langle {H}^{(m-1)}\rangle)_{d-d_m}=0$.
Therefore, the sequence \eqref{eq:H1exact} implies $H_1(K_{\bullet}(\bm{H}))_d=0$.

\medskip\noindent
\emph{Induction step for the second assertion.}
Assume $\max\{d_m, d_{\rm reg}(\langle {H}^{(n+1)}\rangle)\}\le D^{(n;d_1,\ldots,d_{n+1})}$.
We argue again by induction on $m$.
Let $m>n+1$ and $d\ge D^{(n;d_1,\ldots,d_n)}+d_{n+1}$.
Since $d_{m-1}\leq d_m$, the additional hypothesis is inherited by $\bm H^{(m-1)}$:
$\max\{d_{m-1}, d_{\rm reg}(\langle H^{(n+1)}\rangle)\} \leq D^{(n;d_1,\ldots,d_{n+1})}$.
Hence the induction hypothesis gives $H_1(K_\bullet(\bm H^{(m-1)}))_d=0$, so it suffices to show $H_0(K_\bullet(\bm{H}^{(m-1)}))_{d-d_m}=0$.
By $d_m\le D^{(n;d_1,\ldots,d_{n+1})}$, we have
{\small 
\vskip -5mm
\begin{align*}
d-d_{m}& \geq  D^{(n;d_1,\ldots,d_n)}+d_{n+1}-d_{m} \geq D^{(n;d_1,\ldots,d_n)}+d_{n+1}-D^{(n;d_1,\ldots,d_{n+1})} \\
& =  (d_1+\cdots +d_{n+1}-n-1)+2-\left(\left\lfloor \frac{d_1+\cdots +d_{n+1}-n-1}{2}\right\rfloor+1\right)\\
& \geq  \left\lfloor \frac{d_1+\cdots +d_{n+1}-n-1}{2}\right\rfloor+1=D^{(n;d_1,\ldots,d_{n+1})} \geq d_{\rm reg}(\langle H^{(n+1)} \rangle).
\end{align*}}

\vskip -1mm
\noindent
Since $\langle {H}^{(n+1)}\rangle\subset \langle {H}^{(m-1)}\rangle$, we also have $d_{\rm reg}(\langle{H}^{(m-1)}\rangle)\le d_{\rm reg}(\langle {H}^{(n+1)}\rangle)$, and therefore\\
$H_0(K_{\bullet}(\bm{H}^{(m-1)}))_{d-d_m}=(R/\langle {H}^{(m-1)}\rangle)_{d-d_m}=0$, as desired.
\end{proof}

\begin{proposition}\label{prop:sdb}
Let $\bm{F} = (f_1,\ldots,f_m)$ be a sequence of non-constant polynomials in $R$ with $m >  n$ and $1\leq d_{n+1} \leq \cdots \leq d_m$. 
Assume that $(\bm{F}^{\rm top})^{(n)}$ forms a regular sequence in $R$.
Then we have $\widetilde{d}_{\rm reg}(\langle F^h\rangle)\leq D^{(n;d_1,\ldots,d_n)}+d_m-1$.
Moreover, if $\max\{d_m, d_{\rm reg}(\langle (F^{\rm top})^{(n+1)}\rangle)\}\leq D^{(n;d_1,\ldots,d_{n+1})}$, then $\widetilde{d}_{\rm reg}(\langle F^h\rangle) \leq D^{(n;d_1,\ldots,d_n)}+d_{n+1}-1 = D^{(n+1;d_1,\ldots,d_{n+1})}$.
\end{proposition}

\begin{proof}
Let $\bm{F}'=(f_1^h,\dots,f_m^h,y)$.
Regarding the Koszul complex $K_{\bullet}(\bm{F}')$ as the mapping cone of $K_{\bullet}(\bm{F}^h)(-1)\xrightarrow{\times y}K_{\bullet}(\bm{F}^h)$, we obtain a long exact sequence of graded $K$-vector spaces
\begin{equation}
\xymatrix{
\cdots \ar[r] &
H_1(K_{\bullet}(\bm{F}^h))_{d-1}
\ar[r]^{\times y} &
H_1(K_{\bullet}(\bm{F}^h))_{d}
\ar[r] &
H_1(K_{\bullet}(\bm{F}'))_{d}
\ar[dll]  \\
&
H_0(K_{\bullet}(\bm{F}^h))_{d-1}
\ar[r]_{\times y} &
H_0(K_{\bullet}(\bm{F}^h))_{d}
\ar[r] &
H_0(K_{\bullet}(\bm{F}'))_{d}
\ar[r] & 0 .
}
\label{eq:longexactseq}
\end{equation}

We first show that $H_1(K_\bullet(\bm F'))_d=0$ for all $d\ge \hat d$, where $\hat d=D^{(n;d_1,\ldots,d_n)}+d_m$ in the first assertion and $\hat d=D^{(n;d_1,\ldots,d_n)}+d_{n+1}$ in the second one.
Here we put $\bm H=(f_1^h,\ldots,f_n^h,y,f_{n+1}^h,\ldots,f_m^h)$.
Recall that, for each $k$ with $1\leq k \leq m$, the specialization $y\mapsto 0$ induces a natural isomorphism of graded $K$-algebras
\begin{equation}\label{eq:key_isom}
    R'/\langle f_1^h,\ldots,f_k^h,y\rangle \cong R/\langle f_1^{\rm top},\ldots,f_k^{\rm top}\rangle.
\end{equation}
In particular, the corresponding Hilbert series coincide.
Since $(\bm F^{\rm top})^{(n)}$ is a regular sequence, the Hilbert-series characterization of regular sequences \cite[Proposition~1.7.4]{Bar} gives
\[
{\rm HS}_{R'/\langle \bm{H}^{(n+1)}\rangle}
={\rm HS}_{R/\langle (\bm{F}^{\rm top})^{(n)} \rangle}
=\frac{\prod_{j=1}^n (1-z^{d_j})}{(1-z)^n}
=\frac{\prod_{j=1}^n (1-z^{d_j})(1-z)}{(1-z)^{n+1}}, 
\]
where $\bm H^{(n+1)}=(f_1^h,\ldots,f_n^h,y)$. 
Therefore, the sequence $\bm{H}^{(n+1)}$ is also regular in $R'$. 
Moreover, we have $d_{\rm reg}(\langle f_1^h,\ldots,f_n^h,y,f_{n+1}^h\rangle) =d_{\rm reg}(\langle(F^{\rm top})^{(n+1)}\rangle)$.
In addition, it follows immediately from the definition of $D^{(n;d_1,\ldots,d_m)}$ (see \eqref{eq:newbound0}) that $D^{(n+1;d_1,\ldots,d_n,1)}=D^{(n;d_1,\ldots,d_n)}$ and $D^{(n+1;d_1,\ldots,d_n,1,d_{n+1})}=D^{(n;d_1,\ldots,d_{n+1})}$.
Hence Lemma~\ref{lem:key}, applied to the homogeneous sequence $\bm H$ in $R'$ with $m+1 > n+1$, gives $H_1(K_\bullet(\bm H))_d=0$ for all $d\ge \hat d$, in each of the first and second assertions.
Since the Koszul complex is unchanged up to the usual sign isomorphism under a permutation of the sequence, we have $H_1(K_\bullet(\bm F'))_d\cong H_1(K_\bullet(\bm H))_d=0$ for all $d\ge \hat d$.

Next, by~\eqref{eq:key_isom} with $k=m$, one has $H_0(K_\bullet(\bm F'))_d \cong(R/\langle F^{\rm top}\rangle)_d$.
Since $d_{\rm reg}(\langle F^{\rm top}\rangle) \leq d_{\rm reg}(\langle(F^{\rm top})^{(n)}\rangle) = D^{(n;d_1,\ldots,d_n)}$, this vector space vanishes for every $d\geq D^{(n;d_1,\ldots,d_n)}$, and hence for every $d\geq\hat d$.

Since $H_1(K_\bullet(\bm F'))_d=H_0(K_\bullet(\bm F'))_d=0$ for every $d\geq\hat d$, exactness of \eqref{eq:longexactseq} shows that multiplication by $y$ induces an isomorphism
\[
H_0(K_\bullet(\bm F^h))_{d-1}\xrightarrow{\ \sim\ }H_0(K_\bullet(\bm F^h))_d
\]
for every $d\geq\hat d$.
Therefore the Hilbert function of $R'/\langle F^h\rangle$ is constant from degree $\hat d-1$ onward, and hence $\widetilde d_{\rm reg}(\langle F^h\rangle)\leq\hat d-1$.
\end{proof}

We can combine Proposition~\ref{prop:sdb} and Proposition~\ref{prop:new} to obtain the following bound.

\begin{theorem}\label{thm:maxGB}
Let $\bm{F} = (f_1,\ldots,f_m)$ be a sequence of non-constant polynomials in $R$ with $m>n$ and $1\leq d_{n+1} \leq \cdots \leq d_m$. 
Assume that $(\bm{F}^{\rm top})^{(n)}$ forms a regular sequence in $R$. 
Then, for any DRL ordering $\prec$ on $R$ (and its homogenization $\prec^h$ on $R'$), we have
\[
{\rm max.GB.deg}_{\prec^h}(F^h) \leq D^{(n;d_1,\ldots,d_n)}+d_m-1= d_1+\cdots+d_{n}+d_m-n.
\]
Moreover, if $\max\{d_m, d_{\rm reg}(\langle (\bm{F}^{\rm top})^{(n+1)}\rangle)\}\leq D^{(n;d_1,\ldots,d_{n+1})}$, then
\[
{\rm max.GB.deg}_{\prec^h}(F^h) \leq D^{(n;d_1,\ldots,d_n)}+d_{n+1}-1=d_1+\cdots+d_n+d_{n+1}-n.
\]
\end{theorem}

\begin{proof}
Set $\hat{d}:=D^{(n;d_1,\ldots,d_n)}+d_m$ (or $\hat{d}:=D^{(n;d_1,\ldots,d_n)}+d_{n+1}$ under the additional assumption).
By Proposition~\ref{prop:sdb}, we have $\widetilde{d}_{\rm reg}(\langle F^h\rangle)\le \hat{d}-1$.
Moreover, since $(\bm{F}^{\rm top})^{(n)}$ is regular, it follows that $d_{\rm reg}(\langle F^{\rm top}\rangle)\le D^{(n;d_1,\ldots,d_n)}\le \hat{d}-1$.
Hence $\max\{\,d_{\rm reg}(\langle F^{\rm top}\rangle),\ \widetilde{d}_{\rm reg}(\langle F^h\rangle)\,\} \le \hat{d}-1$.
Applying Proposition~\ref{prop:new}~{\rm (2)} with $H=F^h$ and $\bm H=\bm F^h$ yields the stated bound on ${\rm max.GB.deg}_{\prec^h}(F^h)$.
\end{proof}

\begin{remark}
Assume that $K$ is an infinite field and fix positive integers $d_1\leq\cdots\leq d_m$.
Then the regularity assumption in Theorem~\ref{thm:maxGB} is generic with respect to the choice of the highest-degree homogeneous parts.
Indeed, the set of regular sequences is a nonempty Zariski-open
subset $U\subset R_{d_1}\times\cdots\times R_{d_n}$.
Hence, the inverse image of $U$ under the natural projection $\pi:R_{d_1}\times\cdots\times R_{d_m}\longrightarrow R_{d_1}\times\cdots\times R_{d_n}$ is again nonempty Zariski-open.
Thus, for a generic sequence of homogeneous polynomials of the prescribed degrees, the first $n$ elements form a regular sequence.
\end{remark}

\begin{corollary}
Let $\bm F=(f_1,\ldots,f_m)$ be a sequence of non-constant polynomials in $R$ with $m>n$, and assume that $\bm F^{\rm top}$ is semi-regular.
After permuting only the generators $f_{n+1},\ldots,f_m$, we may assume that $d_{n+1}\leq\cdots\leq d_m$.
Under this ordering, the first bound in Theorem~\ref{thm:maxGB} applies.
The second bound applies whenever its additional numerical hypothesis is satisfied.
\end{corollary}

\begin{proof}
The characterization of semi-regularity recalled in Section~\ref{sec:pre} shows that $(f_1^{\rm top},\ldots,f_n^{\rm top})$ is regular.
The result follows from Theorem~\ref{thm:maxGB}.
\end{proof}

\begin{remark}
Suppose, in addition, that the generators are indexed in
nondecreasing degree order, $d_1\leq\cdots\leq d_m$, and that $\bm F^{\rm top}$ is semi-regular in this order.
Then the regular initial subsequence $(f_1^{\rm top},\ldots,f_n^{\rm top})$ consists of the $n$ generators of smallest degrees, and the first bound in Theorem~\ref{thm:maxGB} involves the $n$ smallest input degrees together with the largest one.
This should be compared with Lazard's bound (Theorem~\ref{thm:Lazard_zero_dim}), which is formulated after arranging the input degrees in nonincreasing order.
When $d_1=\cdots=d_m$, the two bounds coincide.

Since semi-regularity need not be preserved under an arbitrary
permutation (see \cite[p.~582]{Pardue}), the nondecreasing degree ordering is an additional assumption on the chosen order of the sequence.
\end{remark}

\subsection{New bound using saturation exponent}\label{subsec:satexp}

In this subsection, we consider a set $F$ of non-constant polynomials in $R$ \textcolor{black}{with $F \neq F^h$}, and derive an upper bound on $\mathrm{max.GB.deg}_{\prec^h}(F^h)$ in terms of the saturation exponent with respect to $y$.
We use the following positive convention:
$s$ denotes the least positive integer such that $\langle F\rangle^h=(\langle F^h\rangle : y^\infty)=(\langle F^h\rangle:y^s)$.
Thus $s=1$ when $\langle F^h\rangle$ is already saturated with respect to $y$.
See \cite[p.\ 81]{Singular} for the definition of saturation exponent. 
Here, for a monomial $t= x_1^{\alpha_1}\cdots x_n^{\alpha_n} y^{k}$ of $R'$, we call $x_1^{\alpha_1} \cdots x_n^{\alpha_n}$ its $X$-part.

\begin{proposition}\label{pro:dsat}
With notation as above, put $D:=d_{\rm reg}(\langle F^{\rm top}\rangle)$.
Then, $\widetilde{d}_{\rm reg}(\langle {F}^h\rangle)
\leq D+s-1$.
Moreover, $\mathrm{max.GB.deg}_{\prec^h}(F^h) \le D+s-1$.
Consequently, the same bound applies to the Macaulay and mutant solving degrees of the homogeneous system $F^h$.
\end{proposition}

\begin{proof}
If $D=\infty$, both assertions are immediate.
Hence assume that $D<\infty$.
Consider the following exact sequence:
\[{\xymatrix{
0 \ar[r] & 
\bigl(R'/\langle F\rangle^h\bigr)(-s)\ar[r]^(0.5){\quad \times y^{s}} &
R'/\langle F^h\rangle \ar[r] & 
R'/(\langle F^h\rangle+\langle y^{s}\rangle) \ar[r] & 
0, }}
\]
where $R'/(\langle F^h\rangle: y^{s})=R'/\langle F\rangle^h$. 
Then, we have 
\[
{\rm HS}_{R'/\langle F^h\rangle}(z)=
{\rm HS}_{R'/(\langle F^h\rangle +\langle y^{s}\rangle)}(z)+z^{s}{\rm HS}_{R'/\langle F\rangle^h}(z).
\]

First, we show ${\rm HF}_{R'/(\langle F^h\rangle +\langle y^{s}\rangle)}(d)=0$ for $d\geq D+s-1$, so that ${\rm HF}_{R'/\langle F^h\rangle}(d)={\rm HF}_{R'/\langle F\rangle^h}(d-s)$.
Let $t$ be any monomial in $R'_d$.
If $y^s \mid t$, then $t\in (\langle F^h \rangle + \langle y^s \rangle)_d$.
Otherwise, the exponent of $y$ in $t$ is at most $s-1$, and hence the degree of the $X$-part of $t$ is at least $d-(s-1)\geq D$.
Since \textcolor{black}{$D<\infty$}, recall from \cite[Lemma 3]{KY} that $\langle \mathrm{LM}_{\prec^h}(\langle F^h\rangle)\rangle$ contains any monomial of $R$ of degree $D$, and thus it also contains $t$.
Hence we have $t \in \langle \mathrm{LM}_{\prec^h}(\langle F^h\rangle +\langle y^{s}\rangle)\rangle$ in either case.
This implies that any element in $R'_d$ is reduced to $0$ modulo a Gr\"{o}bner basis of $\langle F^h \rangle + \langle y^s \rangle$ with respect to $\prec^h$, so that $R'_d/(\langle F^h\rangle +\langle y^{s}\rangle)_d = 0$.

Next we show that ${\rm HF}_{R'/\langle F\rangle^h}(d)$ becomes constant for $d\geq D-1$, which implies that ${\rm HF}_{R'/\langle F^h\rangle}(d)$ becomes constant for $d\geq D+s-1$. 
Let $G$ be the reduced Gr\"obner basis of $\langle F\rangle$ with respect to $\prec$. 
Then $G^h$ is a Gr\"obner basis of $\langle F\rangle^h$ with respect to $\prec^h$. 
By \cite[Remark 15]{CG20} (see also \cite[Remark 9]{KY}), we have ${\rm max.GB.deg}_{\prec}(F)\leq D$, and thus any element of $G^h$ is of total degree not greater than $D$. 
Since $R/\langle F^{\rm top}\rangle$ is Artinian, the affine ideal $\langle F \rangle$ is zero-dimensional (see \cite[Remark 6]{KY}), and hence $R/\langle F \rangle$ is finite-dimensional over $K$.
Let $\{t_1,\ldots,t_r\}$ be the standard monomial basis of $R/\langle F\rangle$ as a $K$-vector space, that is,
\[
\{t_1,\ldots,t_r\}=\{t \in R: \mbox{$t$ is a monomial such that } {\rm LM}_{\prec}(g) \nmid  t \mbox{ for any }g\in G\}
\]
with $r:= \dim_K R/\langle F \rangle < \infty$.
It follows from \cite[Lemma 3]{KY} that any monomial of $R_D$ is divisible by $\mathrm{LM}_{\prec}(g)$ for some $g \in G$.
This implies that the degree of each $t_i$ is smaller than $D$.
By Macaulay's basis theorem (cf.\ \cite[Theorem 1.5.7]{KR1}), as a basis of the $K$-vector space $(R'/\langle F\rangle^h)_d$, we can take
\[
LB_d'=\{t\in R'_d : \mbox{$t$ is a monomial such that } {\rm LM}_{\prec^h}(g) \nmid t \mbox{ for any } g\in G^h\}.
\]
For each $d$ with $d\geq D-1$, it is clear that $LB_{d}' \supset \{t_1y^{k_1},\ldots,t_ry^{k_r}\}$ with $\deg(t_i)+k_i=d$, and we claim that the equality holds.
Indeed, for every $t \in LB_d'$, it follows from $\mathrm{LM}_{\prec}(G) = \mathrm{LM}_{\prec^h}(G^h)$ that $\mathrm{LM}_{\prec}(g)\nmid t$ for any $g \in G$.
Therefore, the $X$-part of $t$ coincides with some $t_i$.
Thus, for $d\geq D-1$, it follows that {$\dim_K (R'/\langle F\rangle^h)_d$} is equal to the constant $r$.

Consequently, the Hilbert function of $R'/\langle F^h\rangle$ is constant from degree $D+s-1$ onward.
Hence $\widetilde d_{\rm reg}(\langle F^h\rangle)\leq D+s-1$, which proves the first assertion.
Since $s\geq1$, one has $D\leq D+s-1$.
Therefore Proposition~\ref{prop:new}~{\rm (2)}, together with the
first inequality, gives $ \mathrm{max.GB.deg}_{\prec^h}(F^h) \leq D+s-1$.
\end{proof}

\begin{example}
In \ref{app:example} below, we consider the case where $n=3$, $m=4$, $d_1=d_2=d_3=d_4=2$ and $R=\mathbb{F}_{73}[x_1,x_2,x_3]$, where $F$ has one affine zero $(0,0,0)$ with multiplicity 1. 
In this example, $D=d_{\rm reg}(\langle F^{\rm top} \rangle)=3$, the saturation exponent $s$ is $3$, $\widetilde{d}_{\rm reg}(\langle F^h \rangle)=4$ and ${\rm max.GB.deg}_{\prec^h}(F^h)=4$.
Thus, as $D+s-1=5$, the gap between $D+s-1$ and $\widetilde{d}_{\rm reg}(\langle F^h\rangle)$ is $1$. 
Their Hilbert series are as follows:
\begin{eqnarray*}
{\rm HS}_{R'/\langle F^h\rangle} & = & 1+4z+6z^2+4z^3+z^4+z^5+\cdots,\\
{\rm HS}_{R'/(\langle F^h\rangle +\langle y^3\rangle)} & = & {1+  4z + 6z^2 + 3z^3},\\
{\rm HS}_{R'/\langle F\rangle^h} & = & 1+z+z^2+\cdots.
\end{eqnarray*}
These Hilbert series also verify the identity ${\rm HS}_{R'/\langle F^h\rangle} ={\rm HS}_{R'/(\langle F^h\rangle +\langle y^3\rangle)}+z^3{\rm HS}_{R'/\langle F\rangle^h}$. 
\end{example}

\section{Correspondence and degree falls in Gr\"{o}bner basis computations}\label{sec:app}

We continue to use the notation introduced in Sections~\ref{sec:pre} and~\ref{sec:3}.
In this section, we compare the Gr\"obner basis computations for $\bm F^{\rm top}$, $\bm F^h$, and $\bm F$.
The basic degree-by-degree relation between these computations is already implicit in complexity analyses based on semi-regularity.
Here we need a more explicit description because the bounds obtained at the end of this section concern the solving degree of the particular Buchberger-like computation considered here and depend on the selection and reduction processes.

We therefore consider the Buchberger-like computations specified below and examine when pair selection, reduction, and insertion of remainders commute with specialization.
The computations are {\em synchronized} before the first {\em degree fall}, while the first degree fall marks the point at which full synchronization may break.
We characterize this point in terms of the injectivity of multiplication by the homogenizing variable $y$.
The three systems are naturally related via specialization: $F$ is obtained from $F^h$ by setting $y=1$, whereas $F^{\rm top}$ is obtained by setting $y=0$.
Examples illustrating the correspondence are given in \ref{app:example}.

Let $G$, $G_{\rm hom}$, and $G_{\rm top}$ be the reduced Gr\"obner bases of $\langle F\rangle$, $\langle F^h\rangle$, and $\langle F^{\rm top}\rangle$, respectively, where the monomial orderings are DRL $\prec$ on $R$ and its homogenization $\prec^h$ on $R'=R[y]$.
It is well-known that a Gr\"obner basis of $\langle F^h\rangle$ yields a Gr\"obner basis of $\langle F\rangle$ after dehomogenization by setting $y=1$.

\paragraph{Computational setting}
Throughout this section, we consider a Buchberger-like Gr\"obner basis computation that repeatedly selects a single S-pair, forms the corresponding S-polynomial, and reduces it by the current intermediate basis;
in particular, we do not consider matrix-based simultaneous reductions such as those used in the $F_4$ algorithm.
When comparing the computations for $F^h$ and $F^{\rm top}$, the input generators are kept indexed by $i$ when S-pairs are formed;
thus generators with different indices are not identified even if $f_i^{\rm top}=f_j^{\rm top}$.

We use the normal selection strategy (see the first paragraph of Section~\ref{subsec:complexity}), except in Section~\ref{subsec:overD}, where the choices of S-pairs and reducers are described separately.
Whenever a nonzero remainder is inserted into the current basis, all new S-pairs involving it are generated.
\textcolor{black}{We may use any criterion for discarding unnecessary S-pairs which can be examined only by leading monomials; for example, Buchberger's product criterion.}
Thus an S-pair whose leading monomials are relatively prime is discarded before its S-polynomial is formed.
A discarded S-pair is not regarded as processed; its degree is not counted as a step degree, and the corresponding multiples do not contribute to the algorithmic solving degree.
We do not prescribe how ties among S-pairs are resolved or in which order reducers are chosen.
Proposition~\ref{prop:sync_partial} shows that, under its hypothesis, corresponding S-pairs can be selected in either computation and that, after a reduction sequence is chosen in either computation, a corresponding sequence can be constructed in the other; when $a=0$, reduction steps that specialize to the identity are omitted.

For either of the pairs $(R,\prec)$ and $(R',\prec^h)$, we use the standard notion of polynomial reduction.
Let $A$ denote either $R$ or $R'$ and $\prec_A$ the corresponding monomial ordering.
A nonzero polynomial $g\in A$ is said to be reducible by a nonzero polynomial $h\in A$ if $g$ has a term $ct$ with monomial $t$ and non-zero constant $c$ such that $t$ is divisible by ${\rm LM}_{\prec_A}(h)$. 
In this case, a one-step reduction of $g$ by $h$ replaces $g$ by $g'=g-c\,\frac{t}{{\rm LT}_{{\prec_A}}(h)}\,h$, by which the term $ct$ is eliminated.
Here we write $g\xlongrightarrow[h]{}g'$.
For a set $H$ of nonzero polynomials in $A$, we say that $g$ is {\em reducible by $H$} if it is reducible by some $h\in H$.
A reduction sequence of $g$ by $H$ means a sequence $g=r_0 \to r_1 \to \cdots \to r_k$ such that each $r_{i+1}$ is obtained from $r_i$ by a one-step reduction by some element of $H$.
If $r_k$ is not reducible by $H$, then we call $r_k$ a {\em remainder} of $g$ (on division) by $H$ and write $g\xlongrightarrow[H]{*}{r_k}$.

We denote by $\mathcal G_{\rm hom}$ and $\mathcal G$ the sets of all polynomials that occur as elements of the current intermediate basis during the computations of $G_{\rm hom}$ and $G$, respectively, and similarly define $\mathcal G_{\rm top}$ for $G_{\rm top}$.
Thus each of these sets consists of the original input generators together with all nonzero polynomials inserted into the current basis during the corresponding computation.
These sets need not be reduced; after termination, their inter-reductions give $G_{\rm hom}$, $G$, and $G_{\rm top}$, respectively.

\paragraph{Leading monomials of intermediate bases}
For homogeneous input, the step degrees are nondecreasing under the normal selection strategy, and every nonzero remainder is reduced by the current basis before being inserted.
Hence the leading monomial of every such remainder is minimal with respect to divisibility when it is inserted and remains minimal thereafter.
Consequently, the only possible nonminimal monomials in $\operatorname{LM}_{\prec^h}(\mathcal G_{\rm hom})$ and $\operatorname{LM}_{\prec}(\mathcal G_{\rm top})$ are those contributed by the original input generators.
After removing from these sets the input leading monomials that are properly divisible by another monomial in the same set, the remaining sets are precisely $\operatorname{LM}_{\prec^h}(G_{\rm hom})$ and $\operatorname{LM}_{\prec}(G_{\rm top})$, respectively.

For the affine computation, the leading monomial of every nonzero remainder likewise remains minimal with respect to divisibility until the first degree fall.
Indeed, during this period every nonzero remainder has total degree equal to its step degree, and hence the step degrees are nondecreasing.

\subsection{Degree falls and $y$-injectivity}\label{subsec:degfall}

We introduce three conditions relating the injectivity of multiplication by $y$ to the absence of $y$-divisible leading monomials in the homogeneous computation.
These conditions will later be used to control degree falls in the affine computation.
Let $\mathcal{B}$ be the minimal generating set of the monomial ideal $\langle \mathrm{LM}_{\prec^h}(\langle F^h \rangle)\rangle$, so that $\mathcal{B} = \mathrm{LM}_{\prec^h}(G_{\rm hom})$.
Consider the following conditions for a fixed positive integer $d$:
\begin{enumerate}
\item[(A)] The following multiplication-by-$y$ map is injective:
    \begin{equation}\label{eq:mult-y}
    (R'/\langle F^h \rangle)_{d-1}\longrightarrow (R'/\langle F^h \rangle)_{d} \ ; \ h+\langle F^h\rangle \longmapsto yh+\langle F^h \rangle,
    \end{equation}
    where $R'=R[y]$.
\item[(B)] No monomial in $\mathcal B_d:=\mathcal{B}\cap R'_d$ is divisible by $y$.
\item[(C)] In a Buchberger-like computation with the normal selection strategy for $\langle F^h\rangle$, every nonzero remainder obtained at step degree $d$ (see the beginning of Section \ref{subsec:complexity}) has leading monomial not divisible by $y$.
\end{enumerate}

\begin{lemma}\label{lem:char-no-y}
For each positive integer $d$, the condition {\rm (A)} implies condition {\rm (B)}, and the conditions {\rm (B)} and {\rm (C)} are equivalent.

Conversely, if the condition {\rm (B)} holds for every degree $d<d_0$, then the condition {\rm (A)} also holds for every degree $d<d_0$.
\end{lemma}

\begin{proof}
    (A) $\Longrightarrow$ (B): 
    Suppose that ${\rm LM}_{\prec^h}(g)$ were divisible by $y$ for some $g \in (G_{\rm hom})_{d}$.
    Since $y$ is the smallest variable for the DRL ordering $\prec^h$, this implies that $g$ itself is divisible by $y$, say $g=yh$ with $h\in R[y]_{d-1}$.
    As $g\in \langle F^h\rangle$, the injectivity condition implies that $h\in \langle F^h\rangle_{d-1}$.
    Hence $\mathrm{LM}_{\prec^h}(h)$ is divisible by $\mathrm{LM}_{\prec^h}(g_1)$ for some $g_1\in (G_{\rm hom})_{\leq d-1}$, and so is $\mathrm{LM}_{\prec^h}(g)=y\,\mathrm{LM}_{\prec^h}(h)$.
    This contradicts the minimality of $\mathcal{B}$.
    
    (B) $\Longleftrightarrow$ (C):
    By the observation on leading monomials above for homogeneous input, the elements of $\mathcal B_d$ not arising from the input generators are precisely the leading monomials of the nonzero remainders obtained at step degree $d$.
    On the other hand, $\operatorname{LM}_{\prec^h}(f_i^h) =\operatorname{LM}_{\prec}(f_i)\in R$ for every input generator $f_i$, so none of the input contributions is divisible by $y$.
    Therefore, {\rm (B)} and {\rm (C)} are equivalent.

    To prove the last statement, assume that $y h \in \langle F^h \rangle_d$ with $h \in R[y]_{d-1}$.
    Let $r$ be a remainder of $h$ by $(G_{\rm hom})_{\leq d-1}$, so that $y r \in \langle F^h \rangle_{d}$.
    Suppose $r\ne0$.
    Then the monomial ${\rm LM}_{\prec^h}(y r) = y\,\mathrm{LM}_{\prec^h}(r)$ belongs to $\langle \mathrm{LM}_{\prec^h}(\langle F^h\rangle)\rangle$, whence it is divisible by some element $t\in\mathcal B$.
    Since the condition {\rm (B)} holds for every degree $<d_0$ and $\deg(t)\le d < d_0$, we have $y\nmid t$.
    Therefore, the monomial $t$ satisfies $t\mid \mathrm{LM}_{\prec^h}(r)$ and $\deg(t)\leq d-1$, contradicting the fact that $r$ is a remainder by $(G_{\rm hom})_{\le d-1}$.
    Thus $r=0$, and hence $h\in\langle F^h\rangle_{d-1}$.    
\end{proof}

For later use, we package the degree-wise conditions above into the following cumulative condition:

\begin{condition}[Injectivity condition below $d_0$]\label{cond:inj}
Let $d_0$ be a positive integer.
For every $d < d_0$, the conditions {\rm (A)}, {\rm (B)}, and {\rm (C)} hold.
\end{condition}

We define the {\it $y$-injectivity threshold} $d_{\rm inj}$ by
\begin{equation}
    d_{\rm inj}:=
\inf \{d \in \mathbb{Z}_{\geq 1} \mid (R'/\langle F^h \rangle)_{d-1} \xrightarrow{\times y} (R'/\langle F^h \rangle)_d \mbox{ is not injective} \}
\label{eq:dinj}
\end{equation}
where $\inf\emptyset=\infty$.
Equivalently, Lemma~\ref{lem:char-no-y} shows that $d_{\rm inj}$ is the supremum of the positive integers $d_0$ for which Condition~\ref{cond:inj} holds.

\paragraph{Relation with saturation and semi-regularity}
Put $I=\langle F^h\rangle$.
For every positive integer $d$, the kernel of \eqref{eq:mult-y} is $((I:y)/I)_{d-1}$.
Thus, if $(I:y)\neq I$, then $d_{\rm inj}-1$ is the least
integer $\ell$ such that $(I:y)_\ell\neq I_\ell$, and otherwise $d_{\rm inj}=\infty$.

A criterion of Bayer--Stillman~\cite[Lemma~1.6]{BS}, in the form recalled in \cite[Lemma~2.1]{Hashemi2012}, shows that a linear form $\ell$ is \emph{generic} for a homogeneous ideal $I$ if and only if $(I:\ell)_d=I_d$ for all sufficiently large $d$.
Moreover, when $\ell$ is generic for $I$, the \emph{satiety} $\operatorname{sat}(I)$ is equivalently the least positive integer $m$ such that $(I:\ell)_d=I_d$ for every $d\ge m$.

If $R/\langle F^{\rm top}\rangle$ is Artinian, then the homogenizing variable $y$ is generic for $I$.
Indeed, the same argument as in the proof of \cite[Theorem~11]{CG20} applies, since $R'/(I+\langle y\rangle) \cong R/\langle F^{\rm top}\rangle$ is Artinian, so there is no projective zero of $I$ on $y=0$.
In particular, if $d_{\rm inj}<\infty$, then $d_{\rm inj}\le \operatorname{sat}(I)$.

The saturation exponent $s$ used in Section~\ref{subsec:satexp} measures a different aspect of the same $y$-torsion phenomenon:
It is the least positive integer such that $y^s\bigl((I:y^\infty)/I\bigr)=0$, or equivalently, such that $y^s g\in I$ for every $g\in(I:y^\infty)$.

\begin{lemma}\label{lem:dinj-bounds}
With notation as above, assume that $R/\langle F^{\rm top}\rangle$ is Artinian, and let $D=d_{\rm reg}(\langle F^{\rm top}\rangle)$.
\begin{enumerate}
\item[(1)] If $d_{\rm inj}<\infty$, then $d_{\rm inj} \le \operatorname{sat}(I) \le \widetilde d_{\rm reg}(I) \le D+s-1$.
\item[(2)] If $\bm F^{\rm top}$ is cryptographic semi-regular, then
$D\le d_{\rm inj}$.
Consequently, if in addition $d_{\rm inj}<\infty$, then $D \le d_{\rm inj} \le \operatorname{sat}(I) \le \widetilde d_{\rm reg}(I) \le D+s-1$.
\end{enumerate}
\end{lemma}

\begin{proof}
    (1) The first inequality follows from the preceding discussion.
    Since $R'/(I+\langle y\rangle)$ is Artinian, one has $\dim(R'/I)\le1$, and hence $\widetilde d_{\rm reg}(I)=\operatorname{hilb}(I)$ by Remark~\ref{rem:hilbreg}.
    Thus it follows from \cite[Lemma~5.6]{Hashemi2010} that $\operatorname{sat}(I)\le \widetilde d_{\rm reg}(I)$.
    The last inequality follows from Proposition~\ref{pro:dsat}.

    (2) If $\bm F^{\rm top}$ is cryptographic semi-regular, then $H_1(K_\bullet(\bm F^{\rm top}))_d=0$ for every $d<D=d_{\rm reg}(\langle F^{\rm top}\rangle)$.
    By Lemma~\ref{lem:koszul-specialization}, the same holds for $H_1(K_\bullet(\bm F'))_d$, where $\bm{F}':=(f_1^h,\ldots,f_m^h,y)$.
    Thus \textcolor{black}{the exact sequence}~\eqref{eq:longexactseq} gives the condition {\rm (A)} for every $d<D$, while Lemma~\ref{lem:char-no-y} gives {\rm (B)} and {\rm (C)}.
    Consequently, $D\leq d_{\rm inj}$.
\end{proof}

The converse of Lemma~\ref{lem:dinj-bounds}~{\rm (2)} does not hold.
Indeed, the condition $D\le d_{\rm inj}$, equivalently Condition~\ref{cond:inj} with $d_0=D$, does not imply cryptographic semi-regularity of $\bm F^{\rm top}$;
see Section \ref{app:example-noncsr} in \ref{app:example}.
Thus Condition~\ref{cond:inj} with $d_0=D$ is strictly weaker than cryptographic semi-regularity.

\paragraph{Relation to degree falls}

We next introduce the notions of a \emph{degree fall} and the \emph{early stage}, which will later be related to the homogeneous conditions above.

\begin{definition}[Degree fall and early stage]\label{def:degfall}
Let $d$ be a positive integer.
In a Buchberger-like Gr\"obner basis computation, we say that a \emph{degree fall occurs at step degree $d$} if a nonzero remainder obtained at step degree $d$ has total degree strictly smaller than $d$.
Under the normal selection strategy, we call the part of the computation preceding the first degree fall the \emph{early stage}.
If no degree fall occurs, the whole computation is regarded as the early stage.
\end{definition}

For homogeneous input, every nonzero remainder obtained at step degree $d$ is homogeneous of degree $d$.
Hence no degree fall occurs, and the whole computation is the early stage.
For arbitrary input, under the normal selection strategy, the step degree is nondecreasing during the early stage.

Lemma~\ref{lem:char-no-y} is formulated purely on the homogeneous side.
Its significance for the specialized computations will become clear in the next subsection, where we show that the computation for $F^h$ can be synchronized with that for $F^h|_{y=a}$ throughout the early stage.
In particular, taking $a=1$ gives the corresponding synchronization between the computations for $F^h$ and $F$.

As will be proved in Corollary~\ref{cor:dinj-fall} below, if $d_{\rm inj}<\infty$, then $d_{\rm inj}$ is precisely the first step degree at which a degree fall occurs in the chosen Buchberger-like computation for $\langle F\rangle$.
If $d_{\rm inj}=\infty$, then no degree fall occurs.
Thus $d_{\rm inj}$ gives a homogeneous characterization of the first degree fall in the affine computation.

\begin{remark}[Comparison with the first fall degree]\label{rem:first-fall}
Following \cite[Remark~1]{KY}, we define the first fall degree in terms of Koszul homology by $d_{\rm ff}
:=
\inf\left\{
d\in\mathbb Z_{\geq0}
\ \middle|\
H_1(K_\bullet(\bm F^{\rm top}))_d\neq0
\right\}$ with $\inf\emptyset=\infty$; cf.~\cite[Definition~1.3]{CG23}.
Then $d_{\rm ff}\leq d_{\rm inj}$.
Indeed, Lemma~\ref{lem:koszul-specialization} gives $H_1(K_\bullet(\bm F'))_d \cong H_1(K_\bullet(\bm F^{\rm top}))_d$, which vanishes for every $d<d_{\rm ff}$, and hence the relevant part of~\eqref{eq:longexactseq} shows that \eqref{eq:mult-y} is injective.
In particular, if $\bm F^{\rm top}$ is cryptographic semi-regular, then $D=d_{\rm reg}(\langle F^{\rm top}\rangle)\leq d_{\rm ff}\leq d_{\rm inj}$.
\end{remark}

\subsection{Synchronization among Gr\"obner basis computations via specialization}\label{subsec:sync-aff}

In this subsection, we establish a general correspondence principle between Buchberger-like Gr\"obner basis computations under specialization.
This principle provides a unified framework to compare the three computations associated with $F^h$, $F$, and $F^{\rm top}$, and will be applied in the subsequent analysis of the three computations.
For $a\in K$, we consider the specialization map $R'=R[y] \longrightarrow R \ ;\ f \longmapsto f^{(a)}:=f|_{y=a}$ and compare Buchberger-like Gr\"obner basis computations for $\langle F^h\rangle \subset R'$ and for its specialization $\langle F^h|_{y=a}\rangle \subset R$ under the normal selection strategy.
The key observation is the following:
As long as all polynomials actually used in the computation---both for forming S-pairs and as reducers---have leading monomials not divisible by $y$, the formation and reduction of S-polynomials are compatible with specialization.

We begin with the following easy lemma:

\begin{lemma}\label{lm:onestepred}
Let $g,h\in R'$ be homogeneous, and assume that $y\nmid \mathrm{LM}_{\prec^h}(h)$.
If $g^{(a)}$ admits a one-step reduction by $h^{(a)}$, then there is a corresponding one-step reduction of $g$ by $h$. 
Conversely, if $g$ admits a one-step reduction by $h$, then there is a corresponding one-step reduction of $g^{(a)}$ by $h^{(a)}$ for $a\not=0$. 
For $a=0$, if the monomial of the term of $g$ eliminated by $h$ is not divisible by $y$, then there is a corresponding one-step reduction of $g^{(a)}$ by $h^{(a)}$, and otherwise, $g^{(a)}\rightarrow g^{(a)}$ corresponds, which may be called a \textcolor{black}{\em trivial reduction}.
\end{lemma}

\begin{proof}
Let
\[
g^{(a)} 
\xlongrightarrow[h^{(a)}]{} r'
:=
g^{(a)}-c\,\frac{u}{\mathrm{LT}_{\prec}(h^{(a)})}\,h^{(a)}
\]
be a one-step reduction in $R$ which eliminates the term $cu$ with a monomial $u$ and a non-zero constant $c$ in $g^{(a)}$. 
Thus, $u$ is divisible by $\mathrm{LM}_{\prec}(h^{(a)})$. 
Since $g$ is homogeneous, it follows that $t=u\,y^{\deg(g)-\deg(u)}$ is a monomial of $g$.
By the $y$-freeness of $\mathrm{LM}_{\prec^h}(h)$, we have $\mathrm{LM}_{\prec^h}(h)=\mathrm{LM}_{\prec}(h^{(a)})$, so that $\mathrm{LM}_{\prec^h}(h)\mid t$.
We note that for $a=0$, $\deg(u)$ should be equal to $\deg(g)$.
Hence $g$ admits the one-step reduction in $R'$;
\[
g \xlongrightarrow[h]{} r:=g-c'\,\frac{t}{\mathrm{LT}_{\prec^h}(h)}\,h,
\]
where we set $c'=c/ a^{\deg(g)-\deg(u)}$ if $a\neq 0$, and $c'=c$ if $a=0$.
In either case, $r$ is a homogeneous polynomial with $r^{(a)}=r'$.

The converse can be shown in the same manner as above: 
Let $g$ admit the following one-step reduction which eliminates the term $ct$ in $g$ with monomial $t$ and coefficient $c$;
\begin{equation}\label{eq:oneredhom}
g \xlongrightarrow[h]{} r:=g-c\,\frac{t}{\mathrm{LT}_{\prec^h}(h)}\,h,
\end{equation}
where $c\in K^\times$ and $t$ is a monomial in $R'$. 
The assumption $y\nmid \mathrm{LM}_{\prec^h}(h)$ implies $\mathrm{LT}_{\prec^h}(h) =\mathrm{LT}_{\prec^h}(h)^{(a)} = \mathrm{LT}_{\prec}(h^{(a)})$, which divides the $X$-part of $t$.
Thus, applying the substitution map at $y=a$ to the sequence \eqref{eq:oneredhom}, we have 
\begin{equation}\label{eq:onereda}
g^{(a)} \xlongrightarrow[h^{(a)}]{} r'
:=g^{(a)}-c\,\frac{t^{(a)}}{\mathrm{LT}_{\prec}{(h^{(a)})}}\,h^{(a)}
\end{equation}
with $\mathrm{LT}_{\prec}(h^{(a)})\mid t^{(a)}$.
When $a\not=0$, the above step \eqref{eq:onereda} is a one-step reduction which eliminates the term $ct^{(a)}$ in $g^{(a)}$. 
When $a=0$, if $t^{(a)}\not=0$, that is, $t$ is not divisible by $y$, the above step \eqref{eq:onereda} is also a one-step reduction. 
In the case where $a=0$ and $t^{(a)}=0$, $g^{(a)}=r'$ holds in the above step \eqref{eq:onereda}.
\end{proof}

\begin{definition}[Current intermediate bases]
For an S-pair $\bm{s}$ appearing in the computation for $\langle F^h\rangle$, we denote by $\mathcal{G}_{\rm hom}^{(\prec \bm{s})}$ and $\mathcal{G}_{\rm hom}^{(\preceq \bm{s})}$ the current intermediate bases just before and just after $\bm{s}$ is processed, respectively.
For an S-pair $\bm{u}$ appearing in the computation for $\langle F^h|_{y=a}\rangle$, the notation $\mathcal{G}_{(y=a)}^{(\prec \bm{u})}$ and $\mathcal{G}_{(y=a)}^{(\preceq \bm{u})}$ is defined in the same way.
\end{definition}

\begin{proposition}\label{prop:sync_partial}
Let $a\in K$.
Suppose that, immediately before the next S-pair is selected, the current intermediate bases in the computations for $\langle F^h\rangle$ and $\langle F^h|_{y=a}\rangle$ correspond under $g\mapsto g^{(a)}:=g|_{y=a}$, and that their S-pairs that have neither been processed nor discarded correspond as well.
For $a=0$, the input generators are distinguished by their indices as specified in the computational setting.
Let $\bm{s}=(g_1,g_2)$ be one of these S-pairs in the homogeneous
computation, and put $\bm{s}^{(a)}=(g_1^{(a)},g_2^{(a)})$.
Assume that every polynomial in $\mathcal{G}_{\rm hom}^{(\prec\bm{s})}$ has leading monomial not divisible by $y$.
Then the following hold.
\begin{enumerate}
    \item[(1)] The least common multiples of the leading monomials of $\bm{s}$ and $\bm{s}^{(a)}$ agree.
    Hence, $\bm{s}$ can be selected under the normal selection strategy if and only if $\bm{s}^{(a)}$ can be selected.
    Thus, after either pair is selected in one computation, the other pair can be selected in the other computation.

    \item[(2)] Suppose that a chosen reduction sequence in the homogeneous computation is $S(g_1,g_2)=r_0 \xlongrightarrow[h_0]{}\cdots \xlongrightarrow[h_{k-1}]{}r_k$, where $r_k$ is a remainder by $\mathcal{G}_{\rm hom}^{(\prec\bm{s})}$.
    After specialization, and after omitting trivial steps when $a=0$,
    this gives a reduction sequence of $S(g_1^{(a)},g_2^{(a)})$ whose remainder is $r_k^{(a)}$.

    \item[(3)]
    Conversely, suppose that a chosen reduction sequence of $S(g_1^{(a)},g_2^{(a)})$ in the specialized computation has remainder $r'$.
    This sequence can be lifted to a homogeneous reduction sequence.
    If $a\neq0$, the lifted sequence ends in a homogeneous remainder $r$ satisfying $r^{(a)}=r'$.
    If $a=0$, it reaches a homogeneous polynomial $\widetilde r$ with $\widetilde r^{(0)}=r'$; further reductions of terms whose monomials are divisible by $y$, if necessary, yield a homogeneous remainder $r$ with $r^{(0)}=r'$.
\end{enumerate}
If the homogeneous remainder is zero, or if its specialization is nonzero, the current intermediate bases after the corresponding pairs are processed satisfy $\mathcal{G}_{(y=a)}^{(\preceq\bm{s}^{(a)})}= \mathcal{G}_{\rm hom}^{(\preceq\bm{s})}|_{y=a}$.
If, in addition, the leading monomial of every new nonzero homogeneous remainder is not divisible by $y$, then the S-pairs that have neither been processed nor discarded also correspond after this step.
\end{proposition}

\begin{proof}
For {\rm (1)}, let $US$ be the set of all S-pairs that have neither been processed nor discarded immediately before $\bm{s}$ is selected in the homogeneous computation; in particular, $\bm{s}\in US$.
By the hypothesis on the current computations, $US^{(a)}:=(\bm{u}^{(a)})_{\bm{u}\in US}$ is precisely the family of S-pairs that have neither been processed nor discarded in the specialized computation, with the input indices retained when $a=0$.
For every S-pair $\bm{u}=(h_1,h_2)\in US$, our assumption gives $y\nmid \mathrm{LM}_{\prec^h}(h_i)$ for $i=1,2$.
Thus we have $\operatorname{LT}_{\prec^h}(h_i)^{(a)}=\operatorname{LT}_{\prec}(h_i^{(a)})$ and $\operatorname{LM}_{\prec^h}(h_i)=\operatorname{LM}_{\prec}(h_i^{(a)})$, and hence
\[
\operatorname{lcm} \bigl(\operatorname{LM}_{\prec^h}(h_1),\operatorname{LM}_{\prec^h}(h_2)\bigr) = \operatorname{lcm}\bigl(\operatorname{LM}_{\prec}(h_1^{(a)}), \operatorname{LM}_{\prec}(h_2^{(a)})\bigr).
\]
Suppose first that $\bm{s}=(g_1,g_2)$ is selected in the homogeneous
computation.
As the normal selection strategy is used, for every
$\bm{u}=(h_1,h_2)\in US$ we have
\[
\operatorname{lcm}
\bigl(\operatorname{LM}_{\prec^h}(g_1),
      \operatorname{LM}_{\prec^h}(g_2)\bigr)
\preceq^h
\operatorname{lcm}
\bigl(\operatorname{LM}_{\prec^h}(h_1),
      \operatorname{LM}_{\prec^h}(h_2)\bigr).
\]
Since all these leading monomials and least common multiples belong to
$R$, the restriction of $\prec^h$ to them agrees with $\prec$.
Therefore, the corresponding inequality holds in the specialized
computation, and $\bm{s}^{(a)}$ is among the S-pairs that may be
selected under the normal selection strategy.
The converse follows by the same argument.
Thus, after either $\bm{s}$ or $\bm{s}^{(a)}$ is selected in one
computation, its counterpart may be selected in the other.
Moreover, $S(g_1,g_2)^{(a)}=S(g_1^{(a)},g_2^{(a)})$.

For {\rm (2)}, assume that the following reduction sequence is chosen:
\begin{equation}\label{eq:remainder-seq}
S(g_1,g_2)=r_0 \xlongrightarrow[h_0]{} \cdots \xlongrightarrow[h_{k-1}]{} r_k
\end{equation}
with remainder $r_k$ by $\mathcal{G}^{(\prec \bm{s})}_{\rm hom}$, where each reducer $h_i\in \mathcal{G}^{(\prec \bm{s})}_{\rm hom}$ is homogeneous and satisfies $y\nmid \mathrm{LM}_{\prec^h}(h_i)$.
Applying Lemma~\ref{lm:onestepred} successively to each step in \eqref{eq:remainder-seq}, specialization at $y=a$ yields \begin{equation}\label{eq:remainder-seq-specified}
S(g_1^{(a)},g_2^{(a)})=r_0^{(a)} \xlongrightarrow[h_0^{(a)}]{} \cdots \xlongrightarrow[h_{k-1}^{(a)}]{} r_k^{(a)}.
\end{equation}
For $a\neq0$, \eqref{eq:remainder-seq-specified} is a reduction sequence.
If $a=0$, removing the trivial reductions from \eqref{eq:remainder-seq-specified} gives a reduction sequence whose final polynomial is still $r_k^{(a)}$.

Next we show that $r_k^{(a)}$ is a remainder by $\mathcal{G}^{(\prec \bm{s}^{(a)})}_{(y=a)}$.
If $r_k^{(a)}$ were reducible by some $h'\in\mathcal{G}^{(\prec \bm{s}^{(a)})}_{(y=a)}$, then the assumed correspondence between the current intermediate bases would give an element $h\in\mathcal{G}^{(\prec \bm{s})}_{\rm hom}$ such that $h^{(a)}=h'$.
By Lemma~\ref{lm:onestepred}, $r_k$ would then be reducible by $h$, contrary to the fact that $r_k$ is a remainder.
Thus $r_k^{(a)}$ is a remainder, proving {\rm (2)}.

\medskip

For {\rm (3)}, let
\begin{equation}\label{eq:redseqata}
S(g_1^{(a)},g_2^{(a)})=r_0' \xrightarrow[h_0^{(a)}]{} \cdots \xrightarrow[h_{\ell-1}^{(a)}]{} r_\ell'
\end{equation}
be a chosen reduction sequence in the specialized computation with remainder $r_\ell'$, where each $h_i$ is a corresponding element of $\mathcal{G}^{(\prec \bm{s})}_{\rm hom}$.
Applying Lemma~\ref{lm:onestepred} successively to the one-step reductions in \eqref{eq:redseqata}, we obtain a homogeneous reduction sequence
\begin{equation}\label{eq:redseqhom}
S(g_1,g_2)= \widetilde r_0 \longrightarrow \widetilde r_1 \longrightarrow \cdots \longrightarrow \widetilde r_\ell,
\end{equation}
where $\widetilde r_i^{(a)}=r_i'$ for every $i$.
If $a\neq0$, in the same manner as in (2), it can be shown that $\widetilde r_\ell$ is a remainder.
If $a=0$, the polynomial $\widetilde r_\ell$ need not be a remainder.
Any term of $\widetilde r_\ell$ that is still reducible must, however, have a monomial divisible by $y$; otherwise the corresponding one-step reduction would remain nontrivial after specialization and would contradict the fact that $r_\ell'$ is a remainder.
Reducing only such $y$-divisible terms eventually yields a homogeneous remainder $r$.
Since each of these additional reductions specializes to the identity at $y=0$, we still have $r^{(0)}= \widetilde r_\ell^{(0)}= r_\ell'$.

\medskip

Finally, if the homogeneous remainder is zero, neither computation inserts a new polynomial.
If its specialization is nonzero, the two computations insert corresponding remainders, and hence their current intermediate bases continue to correspond.
If, moreover, every new nonzero homogeneous remainder has leading monomial not divisible by $y$, then specialization preserves its leading monomial and does not annihilate it.
The newly generated S-pairs therefore correspond, and Buchberger's product criterion gives the same decision for them because relative primality of their leading monomials is preserved.
Together with the assumed correspondence for the previously available S-pairs, this proves the last assertion.
\end{proof}

We now derive degree-wise consequences of Proposition~\ref{prop:sync_partial}.
For $a\in K$ and for $d_0\in\mathbb{Z}_{>0}$, let $\mathcal G_{\rm early, (y=a)}^{(<d_0)}$ (resp.\ $\mathcal G_{\rm hom}^{(<d_0)}$) denote the set of intermediate polynomials appearing in the early stage of the Gr\"{o}bner basis computation for $\langle F^h|_{y=a}\rangle$ (resp.\ $\langle F^h \rangle$) up to step degree less than $d_0$.
Supposing that the condition (B) in Section \ref{subsec:degfall} holds for every degree $d<d_0$, any intermediate polynomial obtained up to step degree less than $d_0$ has no leading monomial divisible by $y$. 
Hence specialization at any $a\in K$, including $a=0$, preserves all intermediate polynomials in the early stage.
This yields the following degree-wise correspondence, where $d_{\rm inj}$ is defined as in \eqref{eq:dinj}:

\begin{corollary}\label{cor:sync_degwise}
Let $d_0$ be a positive integer satisfying Condition~\ref{cond:inj}, equivalently $d_0\leq d_{\rm inj}$, and let $a\in K$.
Then, through all step degrees less than $d_0$, corresponding S-pairs can be selected in the computations for $\langle F^h\rangle$ and $\langle F^h|_{y=a}\rangle$.
After a reduction sequence is chosen in either computation, a corresponding sequence can be constructed in the other as in Proposition~\ref{prop:sync_partial}.
Moreover, for these choices, one has $\mathcal G^{(<d_0)}_{\rm early, (y=a)}=\{\, g^{(a)} : g\in \mathcal G_{\rm hom}^{(<d_0)} \,\}$ and ${\rm LM}_{\prec}(\mathcal G_{\rm early,(y=a)}^{(<d_0)}) ={\rm LM}_{\prec^h}(\mathcal G_{\rm hom}^{(<d_0)})$.
After removing from the two sets of leading monomials above the nonminimal monomials contributed by the input generators, the resulting sets are both precisely $\{\operatorname{LM}_{\prec^h}(g):g\in G_{\rm hom},\ \deg(g)<d_0\}$.
\end{corollary}

\begin{proof}
By Lemma~\ref{lem:char-no-y}, condition {\rm (C)} holds at every degree $d<d_0$.
Starting with the input generators, we apply Proposition~\ref{prop:sync_partial} successively to the S-pairs of step degree less than $d_0$.
Its hypotheses hold initially, with the input indices retained when $a=0$, and its last assertion shows that they remain valid after each step, because every new nonzero homogeneous remainder has leading monomial not divisible by $y$.
Hence corresponding S-pairs and reduction sequences can be chosen through all step degrees less than $d_0$, and the intermediate polynomials correspond under specialization.
Their degrees are also preserved, so no degree fall occurs in this range and the specialized polynomials belong to the early stage.
Since all leading monomials involved are not divisible by $y$, specialization preserves them, which gives the equality of the sets of leading monomials.
The last assertion follows from the observation on leading monomials of intermediate bases at the beginning of this section.
\end{proof}

\begin{remark}\label{rem:deterministic}
The choices in Corollary~\ref{cor:sync_degwise} may also be made deterministically.
Keep the input generators in their index order and order every later basis element by the time at which it is inserted.
Among the S-pairs allowed by the normal selection strategy, choose the pair whose ordered pair of basis indices is lexicographically smallest;
\textcolor{black}{see \cite{sugar} for the practical efficiency of such a strategy}.
In each reduction, eliminate the largest reducible term and, among the reducers whose leading monomial divides that term, use the first one in the current basis.
Assign the same indices to corresponding polynomials under specialization.

With these choices, the two executions in Corollary~\ref{cor:sync_degwise} correspond until the first degree fall, after reductions that specialize to the identity at $y=0$ are omitted.
Indeed, leading monomials and insertion order are preserved.
Moreover, since $y$ is the smallest variable, a term not divisible by $y$ precedes every term of the same total degree that is divisible by $y$.
Thus the same nontrivial term and reducer are selected on both sides;
any remaining reduction of a $y$-divisible term specializes to the identity.
At step degree $D$, the same index order may be used together with the priority prescribed in Definition~\ref{def:r(g1,g2)} below.
\end{remark}

In particular, for $a=0$ and $d_0\le d_{\rm inj}$, we have $\mathcal G_{\rm top}^{(<d_0)} =\mathcal G_{(y=0)}^{(<d_0)} =\mathcal G_{\rm early,(y=0)}^{(<d_0)}$.
We next turn to $a=1$, where the threshold $d_{\rm inj}$ admits the following computational interpretation.

\begin{corollary}\label{cor:dinj-fall}
If $d_{\rm inj}<\infty$, then the first degree fall in the Buchberger-like computation for $\langle F\rangle$ occurs at step degree $d_{\rm inj}$.
If $d_{\rm inj}=\infty$, no degree fall occurs.
Equivalently, for every positive integer $d_0$, Condition~\ref{cond:inj} with this $d_0$ holds if and only if no degree fall occurs at any step degree less than $d_0$ in the affine computation.
\end{corollary}

\begin{proof}
By Corollary~\ref{cor:sync_degwise} and Lemma~\ref{lem:char-no-y}, no degree fall occurs at any step degree $d<d_{\rm inj}$.
Suppose now that $d_{\rm inj}<\infty$.
It also follows from the definition of $d_{\rm inj}$ and Lemma~\ref{lem:char-no-y} that the condition {\rm (C)} first fails at $d=d_{\rm inj}$.
Otherwise condition {\rm (B)} would hold for every $d<d_{\rm inj}+1$, and the converse part of Lemma~\ref{lem:char-no-y} would imply condition {\rm (A)} at $d=d_{\rm inj}$, a contradiction.

Choose the first nonzero homogeneous remainder $r$ at step degree $d_{\rm inj}$ with $y\mid\operatorname{LM}_{\prec^h}(r)$.
All polynomials used before this remainder is produced have leading monomials not divisible by $y$, so Proposition~\ref{prop:sync_partial} applies up to this point.
Thus $r^{(1)}$ is the corresponding nonzero remainder in the computation for $F$.
Since $r$ is homogeneous and $y\mid\operatorname{LM}_{\prec^h}(r)$, one has $y\mid r$, whence $\deg(r^{(1)})<d_{\rm inj}$ with $r^{(1)}\neq 0$.
Hence a degree fall occurs at step degree $d_{\rm inj}$.

If $d_{\rm inj}=\infty$, the first part applies at every finite step degree, as desired.
\end{proof}

Thus the injectivity condition gives the exact condition governing the absence of degree falls below a prescribed threshold.

\subsection{Consequences at step degrees below $D=d_{\rm reg}(\langle F^{\rm top}\rangle)$}

We apply the correspondence principle established in Section \ref{subsec:sync-aff} to analyze the Gr\"obner basis computations at step degrees below $D=d_{\rm reg}(\langle F^{\rm top}\rangle)$ {\em under the assumption that $\bm{F}^{\rm top}$ is cryptographic semi-regular}. 
In these steps, no degree fall occurs and the injectivity condition holds, see Condition \ref{cond:inj}.
We write $\mathcal G_{\rm early}^{(<D)}$ for the sets of intermediate polynomials obtained up to step degree $<D$ in the early stage of $\langle F\rangle$, and also write $\mathcal G_{\rm hom}^{(<D)}$ and $\mathcal G_{\rm top}^{(<D)}$ for the sets of intermediate polynomials obtained up to step degree $<D$ in the computations $\langle F^h\rangle$ and $\langle F^{\rm top}\rangle$, respectively.
When $D=\infty$, the superscript $(<D)$ refers to all intermediate polynomials occurring in the corresponding computation.

\begin{corollary}\label{cor:sync_csr}
Assume that $\bm F^{\rm top}$ is cryptographic semi-regular, and let $D=d_{\rm reg}(\langle F^{\rm top}\rangle)$.
If $D<\infty$, then, before the step degree first reaches $D$, corresponding S-pairs and reduction sequences can be chosen in the three computations for $\langle F^h\rangle$, $\langle F\rangle$, and $\langle F^{\rm top}\rangle$ so that they correspond under specialization, after trivial reductions at $y=0$ are omitted.
If $D=\infty$, the same conclusion holds through every finite step degree.
At each step in this range, every polynomial used for forming S-pairs or as a reducer has leading monomial not divisible by $y$.
Moreover, for these choices, one has $\mathcal G_{\rm early}^{(<D)}=(\mathcal G_{\rm hom}^{(<D)})|_{y=1}$, $\mathcal G_{\rm top}^{(<D)}=(\mathcal G_{\rm hom}^{(<D)})|_{y=0}$, and $\operatorname{LM}_{\prec}(\mathcal G_{\rm early}^{(<D)})=\operatorname{LM}_{\prec^h}(\mathcal G_{\rm hom}^{(<D)})=\operatorname{LM}_{\prec}(\mathcal G_{\rm top}^{(<D)})$.
\end{corollary}

\begin{proof}
The cryptographic semi-regularity of $\bm{F}^{\rm top}$ implies that Condition~\ref{cond:inj} holds in every degree $d<D$.
If $D<\infty$, Corollary~\ref{cor:sync_degwise} applies with $d_0=D$.
If $D=\infty$, it applies with every finite positive integer $d_0$.
Thus the stated choices and the correspondence of intermediate polynomials follow.
\end{proof}

Thus, throughout the range specified in Corollary~\ref{cor:sync_csr}, the affine computation and the computation for $F^{\rm top}$ are obtained from the homogeneous computation by specialization at $y=1$ and $y=0$, respectively.

\begin{proposition}\label{prop:interreduce-top}
With notation as above, $(G_{\rm top})_{<D}=(G_{\rm hom})_{<D}\big|_{y=0}$.
\end{proposition}

\begin{proof}
Since $\bm{F}^{\rm top}$ is cryptographic semi-regular, the condition (C) of Section \ref{subsec:degfall} holds for each step degree $d<D$.
Hence every element of $\mathcal{G}_{\rm hom}^{(<D)}$ is homogeneous and has leading monomial not divisible by $y$.
It also follows from Corollary \ref{cor:sync_csr} that $\mathcal G_{\rm top}^{(<D)}=(\mathcal G_{\rm hom}^{(<D)})|_{y=0}$.
If $D<\infty$, all S-pairs of step degree less than $D$ have been processed; if $D=\infty$, the computations have terminated.
Therefore, the inter-reductions of $\mathcal G_{\rm hom}^{(<D)}$ and $\mathcal G_{\rm top}^{(<D)}$ give $(G_{\rm hom})_{<D}$ and $(G_{\rm top})_{<D}$, respectively.
Let $g\in \mathcal{G}_{\rm hom}^{(<D)}$, and consider a reduction sequence of $g$ by $ \mathcal{G}_{\rm hom}^{(<D)} \smallsetminus \{ g\}$.
By the same argument as in the proof of Proposition~\ref{prop:sync_partial} (using the corresponding reduction sequences under specialization), we obtain the corresponding reduction sequence of $g|_{y=0}$ by $\mathcal{G}_{\rm top}^{(<D)}\smallsetminus\{g|_{y=0}\}$.
Hence the remainder of $g$ by $ \mathcal{G}_{\rm hom}^{(<D)} \smallsetminus\{g\}$ is specialized to the remainder of $g|_{y=0}$ by $\mathcal{G}_{\rm top}^{(<D)}\smallsetminus\{g|_{y=0}\}$.
Applying this argument to all elements, we conclude that the inter-reduction of $\mathcal{G}_{\rm hom}^{(<D)}$ is specialized to that of $\mathcal{G}_{\rm top}^{(<D)}$, so that $(G_{\rm top})_{<D}=(G_{\rm hom})_{<D}\big|_{y=0}$.
\end{proof}

\begin{remark}
The above argument shows that the correspondence between the computations of $\langle F^h\rangle$ and $\langle F^{\rm top}\rangle$ up to degree less than $D$ depends only on the injectivity condition (Condition \ref{cond:inj}).
In particular, the equality ${\rm LM}_{\prec^h}(\mathcal{G}_{\rm hom}^{(<D)})={\rm LM}_{\prec}(\mathcal{G}_{\rm top}^{(<D)})$ implies ${\rm LM}_{\prec^h}(G_{\rm hom})_{<D}={\rm LM}_{\prec}(G_{\rm top})_{<D}$.
This yields an alternative computation-theoretic proof of \cite[Lemma~2]{KY}, based on tracking the computation process.
In contrast to the original proof in \cite{KY}, which relies on Hilbert function arguments and semi-regularity, the above proof requires only Condition~\ref{cond:inj} with $d_0=D$.
\end{remark}


\subsection{Synchronization among Gr\"obner basis computations at step degree $D$}\label{subsec:overD}

Throughout this subsection, we assume that $D=d_{\rm reg}(\langle F^{\rm top}\rangle)<\infty$ and that $\bm F^{\rm top}$ is cryptographic semi-regular.
At the first stage at which the step degree reaches $D$, degree falls may occur, so Proposition~\ref{prop:sync_partial} no longer ensures that an entire affine reduction sequence corresponds to a homogeneous one.
We therefore begin with the affine computation and introduce below the choices of S-pairs and reducers used in this subsection.
We will show that an initial part of each affine reduction lifts to the homogeneous computation.
We denote by $\mathcal G_{\rm hom}^{(D)}$ the set of intermediate polynomials obtained at step degree $D$ in the homogeneous computation.

First we recall the {\em sugar strategy}, which is an extension of the normal selection strategy devised for improving the practicality of the computation, and {\em redefine} the step degree by using the {\em sugar degree}.
This change does not affect $\mathcal{G}_{\rm early}^{(<D)}$, see Remark \ref{rm:falldegree} below.
By this definition, the (redefined) step degrees form a nondecreasing sequence throughout the computation of $\langle F\rangle$, as in the homogeneous computation of $\langle F^h\rangle$.

\begin{definition}[Sugar strategy \cite{sugar}]
For each polynomial $g$ appearing as an {\em intermediate polynomial} in the computation of $\mathcal{G}$, its {\em sugar degree}, denoted by ${\rm sdeg}(g)$ here, is defined in a recursive manner as follows:
\begin{itemize}
    \item[(1)] For each $f_i$, we set ${\rm sdeg}(f_i)=\deg(f_i)$.
    \item[(2)] For a polynomial $g$, a monomial $t$, and a non-zero constant $c$, we set ${\rm sdeg}(ctg)={\rm sdeg}(g)+\deg(t)$.
    \item[(3)] For polynomials $g$ and $h$, we set ${\rm sdeg}(g+h)={\rm max}\{{\rm sdeg}(g),{\rm sdeg}(h)\}$.
\end{itemize}
Thus, for each one-step reduction $r\xlongrightarrow[h]{} r'=r-cth$, where $c\in K^\times$ and $t$ is a monomial, we have ${\rm sdeg}(r')=\max\{{\rm sdeg}(r),{\rm sdeg}(h)+\deg(t)\}$. 
For each S-pair $(g_1,g_2)$ appearing in the computation of $\mathcal{G}$, we define its sugar degree by $\mathrm{sdeg}(S(g_1,g_2))$, namely $\max\{{\rm sdeg}(t_1g_1), {\rm sdeg}(t_2g_2)\}$, where $S(g_1,g_2)=t_1g_1-t_2g_2$ with $t_i=\mathrm{lcm}(\mathrm{LM}_{\prec}(g_1),\mathrm{LM}_{\prec}(g_2))/\mathrm{LT}_{\prec}(g_i)$. 
The {\em sugar strategy} selects S-pairs in order of lowest sugar degree, where the tie-breaking is done by the normal selection strategy. 
A ``sugar degree'' can be considered as a ``phantom degree'', where each polynomial $g$ in $R$ with sugar degree $s$ may be associated to a homogeneous polynomial in $R'$ with degree $s$.
\end{definition}

We use the sugar strategy for the affine computation throughout this subsection.
We denote by $\mathcal G^{(<D)}$ and $\mathcal G^{(D)}$ the sets of intermediate polynomials obtained while S-pairs of sugar degree less than $D$ and equal to $D$ are processed, respectively.
We denote by $\mathcal G_{\rm hom}^{(D)}$ the set of intermediate polynomials obtained at step degree $D$ in the homogeneous computation.
Remark~\ref{rm:falldegree} below shows that $\mathcal G^{(<D)}=\mathcal G_{\rm early}^{(<D)}$.

The following definition specifies the intermediate polynomial retained for use in later reductions.
This convention is related to the \emph{saccharine strategy} introduced in~\cite{sugar} in the context of the \emph{double sugar strategy}, in which every reduction preserves the sugar degree.

\begin{definition}\label{def:r(g1,g2)}
For each S-polynomial $S(g_1,g_2)$, by prioritizing reducers which do not change the sugar degree as far as possible, the following reduction sequence is computed:
\begin{equation}\label{eq:red_in_F}
S(g_1,g_2)=r_0
\underset{h_0}{\longrightarrow} r_1\underset{h_1}{\longrightarrow} \cdots
\underset{h_{k-1}}{\longrightarrow}r_k\underset{h_{k}}{\longrightarrow}r_{k+1}
\longrightarrow \cdots \longrightarrow r_\ell, 
\end{equation}
where ${\rm sdeg}(r_0)=\cdots={\rm sdeg}(r_k)$, and, if $k<\ell$, ${\rm sdeg}(r_k)< {\rm sdeg}(r_{k+1})$.
Note that $r_\ell$ is the computed remainder by the {\em current} basis $\mathcal{G}^{(\prec \bm{s})}:=\mathcal{G}_{(y=1)}^{(\prec \bm{s})}$ with $\bm{s} =(g_1,g_2)$.
Then, we define $r_{\rm sgr}(g_1,g_2)$ by $r_k$ and $r_{\rm fin}(g_1,g_2)$ by $r_\ell$. 
Of course, the sequence \eqref{eq:red_in_F} may stop at $k$, that is, $r_{\rm fin}(g_1,g_2)=r_{\rm sgr}(g_1,g_2)$. 
We may apply the double sugar strategy, where we always stop the reduction at $r_{\rm sgr}(g_1,g_2)$.
\end{definition}

After the reduction of the selected S-polynomial has been completed, each nonzero $r_{\rm sgr}(g_1,g_2)$ is retained in the current basis, as is the nonzero final remainder if it is different from $r_{\rm sgr}(g_1,g_2)$.
The corresponding new S-pairs are generated before the next S-pair is selected.
This rule is used only while S-pairs of sugar degree $D$ are processed.
After all such pairs have been processed, the computation continues from the current basis, and no $r_{\rm sgr}(g_1,g_2)$ distinct from the final remainder is retained.

Using this choice of reducers, we now show that, for each affine S-pair selected at sugar degree $D$, the reduction steps from its S-polynomial to $r_{\rm sgr}(g_1,g_2)$ lift to homogeneous reduction steps.

\begin{remark}\label{rm:falldegree}
For any $g\in \mathcal{G}_{\rm early}^{(<D)}$, it can be shown inductively that ${\rm sdeg}(g)=\deg(g)$ holds by using the correspondence between $\mathcal{G}_{\rm early}^{(<D)}$ and $\mathcal{G}_{\rm hom}^{(<D)}$.
Indeed, each $g\in \mathcal{G}_{\rm early}^{(<D)}$ is not a degree fall element and its counterpart $g^h\in \mathcal{G}_{\rm hom}^{(<D)}$ is not divisible by $y$, see Proposition~\ref{prop:sync_partial}. 
Thus, for each S-pair at step degree less than $D$, its sugar degree coincides with its degree.
Consequently, before the step degree first reaches $D$, the sugar strategy makes the same choices of S-pairs as the normal selection strategy, and every reduction preserves the sugar degree.
In particular, $r_{\rm sgr}=r_{\rm fin}$ at such step degrees.
Also, for our {\em revised} step degree, $\mathcal{G}_{\rm early}^{(<D)}=\mathcal{G}^{(<D)}$ holds.
Moreover, Lemma~\ref{lm:sdeg_of_Spair} shows that every S-pair of sugar degree $D$ has S-pair degree $D$.
\end{remark}

\begin{remark}
The preceding reduction convention is particularly simple for signature-based algorithms~\cite{EF}, including $F_5$-type algorithms.
For a Schreyer ordering (see \cite[p.\ 161]{Singular}) satisfying $t_i\bm e_i\prec t_j\bm e_j$ whenever $\deg(t_if_i)<\deg(t_jf_j)$, where each basis element $\bm{e}_i$ of a free $R$-module of rank $m$ corresponds to $f_i$, the sugar degree agrees with the total degree of the signature, and any one-step reduction does not change such degree. 
Hence $r_{\rm fin}(g_1,g_2)$ always coincides with $r_{\rm sgr}(g_1,g_2)$.
\end{remark}

Now we consider the process of the Gr\"obner basis computation for $F$ at step degree $D$ with the sugar strategy.
In this process, a pair $(g_1,g_2)$ is chosen such that $g_1$ and $g_2$ are already computed {\em intermediate polynomials} and ${\rm sdeg}(S(g_1,g_2))=D$.  
Writing $S(g_1,g_2)=c_1t_1g_1-c_2t_2g_2$ for some monomials $t_1,t_2$ and non-zero constants $c_1,c_2$, we have $\max\{{\rm sdeg}(t_1g_1),{\rm sdeg}(t_2g_2)\}=D$.
Here $\mathrm{sdeg}(g_i) \leq D$ for $i=1,2$; otherwise $\max\{{\rm sdeg}(t_1g_1),{\rm sdeg}(t_2g_2)\}>D$.
Put $F_D:=\{f_i\in F\mid \deg(f_i)=D\}$.

\begin{lemma}\label{lm:sdeg_of_Spair}
For each S-pair $(g_1,g_2)$ of sugar degree $D$, we have ${\rm sdeg}(g_i)=\deg(g_i)<D$ or $g_i\in F_D$ for $i=1,2$, so that the degree of $(g_1,g_2)$ is also $D$.
\end{lemma}

\begin{proof}
Without loss of generality, we may assume ${\rm sdeg}(g_1)\geq {\rm sdeg}(g_2)$. 
If ${\rm sdeg}(g_1)<D$, then ${\rm sdeg}(g_2)<D$ and ${\rm sdeg}(g_i)=\deg(g_i)$ by Remark \ref{rm:falldegree}.
Suppose, to the contrary, that ${\rm sdeg}(g_1)=D$ holds for some $g_1\not\in F$.
(In this case, $t_1=1$ and $g_1$ is computed from some S-polynomial.) 
If ${\rm sdeg}(g_2)<D$, then $g_2$ should belong to $\mathcal{G}_{\rm early}^{(<D)}$ and is computed before $g_1$.
Otherwise, that is, ${\rm sdeg}(g_2)=D$, we may set the pair $(g_1,g_2)$ so that $g_2$ is computed before $g_1$.
(We note that all polynomials of degree $D$ in $F$ are always computed before $g_1$.) 
In either case, we may assume that $g_2$ is computed before $g_1$. 
Then $g_1$ is reducible by $g_2$, and a one-step reduction gives a nonzero scalar multiple of $S(g_1,g_2)$, still of sugar degree $D$.
Since $g_2$ had already been computed when $g_1$ was obtained, this is impossible: $r_{\rm fin}$ is reduced, while $r_{\rm sgr}$ admits no further reduction preserving the sugar degree.
\end{proof}

Conversely, if the degree of S-pair $(g_1,g_2)$ with $g_i \in \mathcal{G}^{(<D)}\cup F_D$ is $D$, then clearly the sugar degree of $(g_1,g_2)$ is also $D$.

Therefore, it follows from Lemma \ref{lm:sdeg_of_Spair} that
\begin{align*}
    \mathsf{Spairs}_D :=& \{(g_1,g_2)\mid {\rm sdeg}(S(g_1,g_2))=D,\ g_1,g_2\in \mathcal{G}^{(<D)} \cup F_D\}\\
    =& \{(g_1,g_2)\mid \mbox{the degree of $(g_1,g_2)$ is $D$},
    \ g_1,g_2\in \mathcal{G}^{(<D)} \cup F_D\}.
\end{align*}
By Corollary \ref{cor:sync_csr}, we have $(\mathcal{G}^{(<D)})^h=\mathcal{G}_{\rm hom}^{(<D)}$, and hence
\begin{align*}
    \mathsf{Spairs}_D^h :=&\{(g_1^{h},g_2^{h})\mid (g_1,g_2)\in \mathsf{Spairs}_D\}\\
    =&\{(g'_1,g'_2)\mid \deg(S(g'_1,g'_2))=D,\ g'_1,g'_2\in \mathcal{G}_{\rm hom}^{(<D)}\cup F_D^h\}.
\end{align*}
When a pair $\bm{s}=(g_1,g_2)\in\mathsf{Spairs}_D$ is selected in the affine computation, we process the corresponding pair $\bm{s}^h:=(g_1^h,g_2^h)\in\mathsf{Spairs}_D^h$ in the homogeneous computation.

\begin{remark} 
In the computation of $\mathcal{G}_{\rm hom}$, for each S-pair $(g'_1,g'_2)$ chosen at the step degree $D$, we have $\deg(g_i')<D$ or $g_i'\in F^h_D$ for $i=1,2$. 
This can be shown in the same manner as Lemma \ref{lm:sdeg_of_Spair}. 
\end{remark}

\begin{lemma}\label{lem:cor-atD}
With notation as above, we have the following:
\begin{itemize}
\item[(1)] For $i<k$, to each one-step reduction 
\begin{equation}\label{eq:onered_in_F}
r_i\xlongrightarrow[h_i]{} r_{i+1}:=r_i-c_it_ih_i
\end{equation}
in the sequence \eqref{eq:red_in_F}, where $c_i\in K^\times$ and $t_i$ is a monomial in $R$, the following one-step reduction corresponds:
\begin{equation}\label{eq:onered_in_F^h}
r_i^*\xlongrightarrow[h_i^*]{} r_{i+1}^*=r_{i}^*-c_it_ih_i^*y^{\delta_i}, 
\end{equation} 
where $r_i^*=y^{a_i}r_i^h$ for $a_i=D-\deg(r_i)$, $r_{i+1}^*=y^{a_{i+1}}r_{i+1}^h$ for $a_{i+1}=D-\deg(r_{i+1})$, $h_i^*=y^{b_i}h_i^h$ for $b_i={\rm sdeg}(h_i)-\deg(h_i)$ and $\delta_i=D-\deg(t_i)-{\rm sdeg}(h_i)$. 
Note that, as the reduction \eqref{eq:onered_in_F} does not change the sugar degree for $i<k$, we have $\delta_i=D-\deg(t_i)-{\rm sdeg}(h_i)\geq 0$.

\item[(2)] The part of the affine reduction sequence~\eqref{eq:red_in_F} from $S(g_1,g_2)$ to $r_k=r_{\rm sgr}(g_1,g_2)$ lifts to the homogeneous reduction sequence
\begin{equation*}
{S(g^h_1,g^h_2)=r_0^*
\xlongrightarrow[h^*_0]{} r^*_1\xlongrightarrow[h^*_1]{} \cdots 
\xlongrightarrow[h^*_{k-1}]{} r^*_k}.
\end{equation*}
Then $r_k^*=r_{\rm sgr}(g_1,g_2)^*$ is the computed remainder by the {\em current} basis $\mathcal{G}_{\rm hom}^{(\prec \bm{s}^h)}$, where $\bm{s}=(g_1,g_2)$.
In particular, $h_0^*,\ldots,h_{k-1}^*$ belong to $\mathcal{G}^{(\prec \bm{s}^h)}_{\rm hom}$. 
\end{itemize}
Therefore, together with the obvious correspondence $f_i^h\leftrightarrow f_i$ for $f_i\in F_D$, (1) and (2) give a one-to-one correspondence between $\mathcal G_{\rm hom}^{(D)}$ and $\widehat{\mathcal{G}}^{(D)}:=\{g\in \mathcal{G}^{(D)} : {\rm sdeg}(g)=D\}$, the set of all intermediate polynomials in $\mathcal{G}^{(D)}$ with sugar degree $D$.
\end{lemma}

\begin{proof}
(1) First we consider $r_i^h$ and $c_it_ih_i^h$, which are the homogenizations of $r_{i}$ and $c_it_ih_i$, respectively. 
Then, by adjusting the degrees by powers of $y$, we have $r_{i+1}^h y^{\epsilon_i}= r_i^h - c_i t_i y^{\deg(r_i)-\deg(t_i)-\deg(h_i)} h_i^h$ for $\epsilon_i=\deg(r_{i})-\deg(r_{i+1})$.
Multiplying both sides by $y^{D-\deg(r_i)}$, we have
\begin{eqnarray*}
r_{i+1}^*(=y^{D-\deg(r_{i+1})}r_{i+1}^h)  & = & y^{D-\deg(r_i)+\epsilon_i}r_{i+1}^{h} \\ 
& = & r_i^\ast - c_i t_i h_i^h y^{D-\deg(t_i)-\deg(h_i)}\\
& = & r_i^\ast - c_i t_i h_i^* y^{D-\deg(t_i)-\deg(h_i)-({\rm sdeg}(h_i)-\deg(h_i))}\\
& = & r_i^\ast - c_i t_i h_i^* y^{D-\deg(t_i)-{\rm sdeg}(h_i)}.\end{eqnarray*}
The one-step reduction \eqref{eq:onered_in_F} eliminates the term $c_i t_i\, {\rm LT}_\prec(h_i)$ of $r_i$, and $r_{i}^{\ast}$ has a term\\
$c_i t_i\, \mathrm{LT}_{\prec}(h_i) y^{D-\deg(t_i)-\deg(h_i)}$.
As $\mathrm{LT}_{\prec^h}(h_i^\ast) = \mathrm{LT}_{\prec}(h_i) y^{\mathrm{sdeg}(h_i)-\deg(h_i)}$, we have\\
$c_i t_i\, \mathrm{LT}_{\prec}(h_i) y^{D-\deg(t_i)-\deg(h_i)} = c_i t_i\, \mathrm{LT}_{\prec}(h_i^\ast) y^{D-\deg(t_i)-\mathrm{sdeg}(h_i)}$.
This implies that \eqref{eq:onered_in_F^h} gives a one-step reduction of $r_i^{\ast}$ by $h_i^{\ast}$, as desired.

\medskip
\noindent
(2) We prove the assertion by induction on the selection order in $\mathsf{Spairs}_D$. 
First we show that every $h_i^*$ belongs to $\mathcal{G}_{\rm hom}^{(\prec \bm{s}^h)}$. 

If ${\rm sdeg}(h_i)<D$, then $h_i$ belongs to $\mathcal{G}^{(<D)}$ ($=\mathcal{G}_{\rm early}^{(<D)}$). 
By the correspondence between $\mathcal{G}^{(<D)}$ and $\mathcal{G}_{\rm hom}^{(<D)}$, we have ${\rm sdeg}(h_i)=\deg(h_i)$ and $h_i^*=h_i^h\in \mathcal{G}_{\rm hom}^{(<D)} \subset \mathcal{G}^{(\prec \bm{s}^h)}_{\rm hom}$ (see Remark \ref{rm:falldegree}). 

If ${\rm sdeg}(h_i)=D$, then $h_i$ belongs to $F$, or we can write $h_i=r_{\rm sgr}(g_3,g_4)$, which is obtained from $S(g_3,g_4)$ for some $(g_3,g_4)$ selected beforehand from $\mathsf{Spairs}_D$.
If $h_i\in F$, it is obvious that $h_i^*=h_i^h\in \mathcal{G}_{\rm hom}^{(\prec \bm{s}^h)}$. 
Thus, suppose that $h_i$ is obtained from $S(g_3,g_4)$.
Since the lemma holds for $S(g_3,g_4)$ by the induction hypothesis, the affine reduction sequence from $S(g_3,g_4)$ to $h_i=r_{\rm sgr}(g_3,g_4)$ lifts to the corresponding homogeneous reduction sequence ending in $h_i^*$.
Thus, $h_i^*$ belongs to {$\mathcal{G}^{(\prec \bm{s}^h)}_{\rm hom}$}. 

Next we show that the reduction sequence from $S(g_1^h,g_2^h)$ stops at $r^*_k$. 
Assume, to the contrary, that $r^*_k$ is reduced to $\hat{r}$ by some element $\hat{h}$ in $\mathcal{G}_{\rm hom}^{(\prec\bm{s}^h)}$.
That is, we have a one-step reduction
\begin{equation}\label{eq:hatr}
\hat{r}=r^*_k-ct\hat{h}
\end{equation}
for some monomial $t\in R'$ and $c\in K^\times$ and $\deg(t)+\deg(\hat{h})= \deg(r^*_k)=D$. 
Then the specialization of \eqref{eq:hatr} gives the following one-step reduction:
\begin{equation}\label{eq:hatr2}
\hat{r}^{(1)}=r_k-ct^{(1)}\hat{h}^{(1)}.
\end{equation}
As $\hat{h}$ belongs to $\mathcal{G}_{\rm hom}^{(\prec \bm{s}^h)}$, the polynomial $\hat{h}^{(1)}$ belongs to $\mathcal{G}^{(<D)}\cup \widehat{\mathcal{G}}^{(D)}$ and ${\rm sdeg}(\hat{h}^{(1)})=\deg(\hat{h})$ by the correspondence between $\mathcal{G}_{\rm early}^{(<D)}$ and $\mathcal{G}_{\rm hom}^{(<D)}$ and by the induction argument on the correspondence between $\widehat{\mathcal{G}}^{(D)}$ and $\mathcal{G}_{\rm hom}^{(D)}$.
Then, 
\[
{\rm sdeg}(ct^{(1)}\hat{h}^{(1)})= \deg(t^{(1)})+{\rm sdeg}(\hat{h}^{(1)})\leq \deg(t)+\deg(\hat{h}) 
= D,
\]
which implies that the reduction \eqref{eq:hatr2} keeps the sugar degree $D$. 
This contradicts the definition of $r_k=r_{\rm sgr}(g_1,g_2)$. 
\end{proof}

\paragraph{Consequence at (redefined) step degree $D$}
As a natural consequence of the correspondence in Lemma \ref{lem:cor-atD}, we have the following lemma. 

\begin{lemma}\label{lem:atD}
The monomial ideal $\langle {\rm LM}_{\prec}(\mathcal{G}^{(\leq D)})\rangle$ contains all monomials of $R$ of degree at least $D$, where $\mathcal{G}^{(\leq D)}:=\mathcal{G}^{(<D)}\cup \mathcal{G}^{(D)}$.
\end{lemma}

\begin{proof}
It suffices to show that every monomial in $R_D$ belongs to $\langle {\rm LM}_{\prec}(\mathcal{G}^{(\leq D)})\rangle$.
Let $t$ be a monomial in $R_D$.
Then, it follows from \cite[Lemma 3]{KY} that $\langle {\rm LM}_{\prec^h}((G_{\rm hom})_{\leq D})\rangle_{R'}\cap R_D=R_D$, whence $t'\mid t$ for some monomial $t'\in  {\rm LM}_{\prec^h}((G_{\rm hom})_{\leq D})$.
If $\deg(t')<D$, the monomial $t'$ belongs to ${\rm LM}_{\prec^h}((G_{\rm hom})_{< D})$, where ${\rm LM}_{\prec^h}((G_{\rm hom})_{< D}) \subset {\rm LM}_{\prec^h}(\mathcal{G}_{\rm hom}^{(< D)})={\rm LM}_{\prec}(\mathcal{G}^{(<D)})$ by Corollary \ref{cor:sync_csr}.
Otherwise, namely if $\deg(t')=D$, then $t=t'\in {\rm LM}_{\prec^h}(G_{\rm hom})_D \subset {\rm LM}_{\prec^h}(\mathcal{G}_{\rm hom}^{(D)})$.

It remains to show ${\rm LM}_{\prec^h}(\mathcal{G}_{\rm hom}^{(D)})\cap R_D \subset \langle {\rm LM}_{\prec}(\mathcal{G}^{(\leq D)})\rangle$.
Let $q\in\mathcal{G}_{\rm hom}^{(D)}$ satisfy ${\rm LM}_{\prec^h}(q)\in R_D$.
If $q\in F_D^h$, then $q=f_i^h$ for some $f_i\in F_D$.
Since ${\rm LM}_{\prec^h}(f_i^h)={\rm LM}_{\prec}(f_i)$ and $f_i\in\mathcal G^{(D)}$, the assertion follows.
Otherwise, by Lemma~\ref{lem:cor-atD}, we have $q=r_{\rm sgr}(g_1,g_2)^*$ for some $(g_1,g_2)\in\mathsf{Spairs}_D$.
Since
\[
r_{\rm sgr}(g_1,g_2)^*
=
y^{D-\deg(r_{\rm sgr}(g_1,g_2))}
r_{\rm sgr}(g_1,g_2)^h
\]
and ${\rm LM}_{\prec^h}(q)$ is not divisible by $y$, we must have $\deg(r_{\rm sgr}(g_1,g_2))=D$.
Hence $q=r_{\rm sgr}(g_1,g_2)^h$ and ${\rm LM}_{\prec^h}(q) ={\rm LM}_{\prec}(r_{\rm sgr}(g_1,g_2))$, where $r_{\rm sgr}(g_1,g_2)\in\mathcal G^{(D)}$.
This proves the desired inclusion.
\end{proof}

\begin{remark}\label{rem:threshold}
The construction above is not tied to the degree of regularity.
More generally, let $T$ be a positive integer such that Condition~\ref{cond:inj} holds in every degree $d<T$.
Suppose that, when the sugar degree first reaches $T$, the reduction convention of Definition~\ref{def:r(g1,g2)} is used while the S-pairs of sugar degree $T$ are processed.
Then the proofs of Lemmas~\ref{lm:sdeg_of_Spair} and \ref{lem:cor-atD} remain valid with $T$ in place of $D$.
Indeed, below $T$ no degree fall occurs, so the sugar degree agrees with the ordinary degree, the sugar strategy selects the same pairs as the normal selection strategy, and every reduction preserves the sugar degree.
After all S-pairs of sugar degree $T$ have been processed, the computation may be continued by Buchberger's algorithm, with only final remainders inserted into the current basis.

In particular, if $d_{\rm inj}<\infty$, one may take $T=d_{\rm inj}$.
If $d_{\rm inj}=\infty$, the preceding agreement holds at every finite step degree and no additional reduction convention is needed.
If, moreover, $\bm F^{\rm top}$ is cryptographic semi-regular, then $D\leq d_{\rm inj}$ and $\mathcal G^{(\leq D)}\subseteq \mathcal G^{(\leq d_{\rm inj})}$.
Here $\mathcal G^{(\leq T)}$ refers to the computation in which the convention of Definition~\ref{def:r(g1,g2)} is applied while S-pairs of sugar degree $T$ are processed.
Consequently, Lemma~\ref{lem:atD} implies that $\bigl\langle\operatorname{LM}_{\prec}(\mathcal G^{(\leq d_{\rm inj})}) \bigr\rangle$ contains every monomial of $R$ of degree at least $d_{\rm inj}$.
\end{remark}

\paragraph{Solving degree of $F$ as the highest S-pair degree}

Finally, Lemma~\ref{lem:atD} allows us to apply the bound of \cite[Lemma~4]{KY} on the highest S-pair degree (equivalently, the highest original step degree) directly to the present computation, without restarting it.

\begin{lemma}[cf.\ {\cite[Lemma 4]{KY}}]\label{lem:sd}
With notation as above, assume that $D\geq {\rm max}\{\deg(f) : f\in F\}$ and that $\bm{F}^{\rm top}$ is cryptographic semi-regular.
Let $\mathcal A$ be a Buchberger-like algorithm satisfying the computational setting above, including Buchberger's product criterion, the sugar strategy, and the reduction convention of Definition~\ref{def:r(g1,g2)}.
Then ${\rm sd}_\prec^{\mathcal{A}} (F) \leq 2D-1$.
If the sugar strategy is stopped after sugar degree $D$ and $\mathcal G^{(\leq D)}$ is inter-reduced before the computation is continued, then the solving degree in the strict sense is at most $2D-2$.
\end{lemma}

\begin{proof}
After all S-pairs of sugar degree $D$ have been processed, Lemma~\ref{lem:atD} shows that the current set $H:=\mathcal G^{(\leq D)}$ satisfies $\langle\mathrm{LM}_{\prec}(H)\rangle_D=R_D$.
This is precisely the property used in the proof of \cite[Lemma~4]{KY}.
Hence that proof applies to the present computation without restarting it and yields ${\rm sd}_{\prec}^{\mathcal A}(F)\le2D-1$.
After the stated inter-reduction of $H$, the second part of the same argument gives the strict bound $2D-2$.
\end{proof}

\subsection*{Acknowledgement}
Many of the contents of this paper were presented at Effective Methods in Algebraic Geometry (MEGA 2024).
The authors thank the anonymous referees assigned by the organization of this conference for their comments and suggestions.
The authors are grateful to Yuta Kambe and Shuhei Nakamura for helpful comments. They are also thankful to Masayuki Noro for providing 
computational tools for our computational experiments. 
This work was supported by JSPS Grant-in-Aid for Young Scientists 20K14301 and 23K12949, JSPS Grant-in-Aid for Scientific Research (C) 21K03377 and 26K06753, JSPS Grant-in-Aid for Scientific Research (B) 25K00140, and JST CREST Grant Number JPMJCR2113, and JST K Program Grant Number JPMJKP24U2, Japan.

\renewcommand{\baselinestretch}{0.85}

\renewcommand{\baselinestretch}{1} 
\appendix
\def\thesection{Appendix \Alph{section}}

\def\thesection{Appendix \Alph{section}}
\section{Examples for Section 4}\label{app:example}

\def\thesection{\Alph{section}}\label{app:A}

We use the same notation as in Section~\ref{sec:app}.
We give two examples illustrating the results of that section.
The first revisits \cite[Example~1]{KY} and illustrates the correspondence among the three Gr\"obner basis computations in the cryptographic semi-regular case, including the behavior at step degree $D$.
The second shows that Condition~\ref{cond:inj} with $d_0=D$ may hold even when $\bm F^{\rm top}$ is not cryptographic semi-regular.

The reduced Gr\"obner bases in the first example were computed with Risa/Asir~\cite{Asir}.
For the affine computation below, we also give a direct Buchberger computation in order to make the intermediate steps explicit.

\subsection{A cryptographic semi-regular example}\label{app:example-csr}

Let $p=73$, $K={\mathbb F}_p$, and
\begin{eqnarray*}
f_1 & = & x_1^2+3x_1 x_2+x_2^2-2x_1x_3-2x_2x_3+x_3^2-x_1-2x_2+x_3,\\
f_2 & = & 4x_1^2+3x_1x_2+4x_1x_3+x_3^2-2x_1-x_2+2x_3,\\
f_3 & = & 3x_1^2+9x_2^2-6x_2x_3+x_3^2-x_1+x_2-x_3,\\
f_4 & = & x_1^2-6x_1x_2+9x_2^2+2x_1x_3-6x_2x_3+2x_3^2-2x_1+x_2
\end{eqnarray*}
with $n=3$ and $m=4$.
Then, $d_1=d_2=d_3=d_4=2$.
The highest-degree homogeneous parts are
{\small 
\begin{eqnarray*}
f_1^{\rm top} & = & x_1^2+3x_1x_2+x_2^2-2x_1x_3-2x_2x_3+x_3^2,\\
f_2^{\rm top} & = & 4x_1^2+3x_1x_2+4x_1x_3+x_3^2,\\
f_3^{\rm top} & = & 3x_1^2+9x_2^2-6x_2x_3+x_3^2,\\
f_4^{\rm top} & = & x_1^2-6x_1x_2+9x_2^2+2x_1x_3-6x_3x_2+2x_3^2.
\end{eqnarray*}}

\noindent
One verifies that $\bm{F}^{\rm top}$ is semi-regular, and hence cryptographic semi-regular.
Moreover, $D=d_{\rm reg}(\langle F^{\rm top}\rangle) =3$.
For ease of comparison with the stepwise computation below, we list the elements of the Gr\"obner bases in increasing order of total degree.
Indeed, the reduced Gr\"obner basis $G_{\rm top}$ of the ideal $\langle F^{\rm top}\rangle$ with respect to the DRL ordering $x_1\succ x_2\succ x_3$ is 
{\small
\[
\{ \underline{x_1^2}+68x_2x_3+55x_3^2,\ 
\underline{x_1x_2}+27x_2x_3+29x_3^2,\ 
\underline{x_2^2}+x_2x_3+71x_3^2,\ 
\underline{x_1x_3}+3x_2x_3+33x_3^2,\ 
\underline{x_2x_3^2},\
\underline{x_3^3}
\}
\]}

\noindent
Therefore, we have ${\rm LM}_{\prec}(G_{\rm top})=\{x_1^2,x_1x_2,x_2^2,x_1x_3,x_2x_3^2,x_3^3\}$, from which the Hilbert-Poincar\'e series of $R/\langle F^{\rm top}\rangle$ is computed as
\[
{\rm HS}_{R/\langle F^{\rm top}\rangle}(z) =1 + 3z + 2z^2 =\left( \frac{(1-z^2)^4}{(1-z)^3}\bmod{z^3} \right),
\]
whence $D=3$. 

On the other hand, the reduced Gr\"obner basis $G_{\rm hom}$ of the ideal $\langle F^h\rangle$ with respect to the DRL ordering $x_1\succ x_2\succ x_3\succ y$ is 
{\small
\begin{eqnarray*}
&& \{
\underline{x_1^2}+68x_2x_3+55x_3^2+72x_1y+40x_2y+14x_3y,\\
&& \underline{x_1x_2}+27x_2x_3+29x_3^2+20x_1y+37x_2y+12x_3y,\\
&& \underline{x_2^2}+x_2x_3+71x_3^2+57x_1y+3x_2y+52x_3y,\\
&& \underline{x_1x_3}+3x_2x_3+33x_3^2+22x_1y+5x_2y+14x_3y,\\
&& \underline{x_2x_3^2}+60x_1y^2+22x_2y^2+39x_3y^2,\quad \underline{x_3^3}+72x_1y^2+14x_2y^2+56x_3y^2,\\
&& \underline{x_2x_3y}+16x_1y^2+55x_2y^2+38x_3y^2,\quad \underline{x_3^2y}+72x_1y^2+66x_2y^2+70x_3y^2,\\
&& \underline{x_1y^3},\ \underline{x_2y^3},\ \underline{x_3y^3}
\},
\end{eqnarray*}}

\noindent
and thus 
\[
{\rm LM}_{\prec^h}(G_{\rm hom})=\{x_1^2,x_1x_2,x_2^2,x_1x_3,x_2x_3^2,x_3^3,x_2x_3y,x_3^2y,x_1y^3,x_2y^3,x_3y^3\}.
\]
Then the Hilbert-Poincar\'e series of $R'/\langle F^h\rangle$ satisfies
\[
\left( {\rm HS}_{R'/\langle F^{h}\rangle}(z) \bmod{z^3} \right)
=\left( 1+4z+6z^2 \bmod{z^3} \right)= \left( \frac{(1-z^2)^4}{(1-z)^4}\bmod{z^3} \right).
\]
We note that ${\rm HF}_{R'/\langle F^h\rangle}(3)=4$ and ${\rm HF}_{R'/\langle F^h\rangle}(4)=1$.
We can also examine the {\em correspondence} ${\rm LM}_{\prec^h}(G_{\rm hom})_{<D}={\rm LM}_{\prec}(G_{\rm top})_{<D}$.
Moreover, no element of ${\rm LM}_{\prec^h}(G_{\rm hom})$ of degree less than $D$ is divisible by $y$.
Thus Condition~\ref{cond:inj} holds with $d_0=D$, and Corollary~\ref{cor:dinj-fall} shows that no degree fall occurs at any step degree less than $D$.

\medskip

We next examine the affine computation more closely, step degree by
step degree.
As in Section~\ref{subsec:overD}, we use the sugar degree as the
step degree in the following computation.
No inter-reduction is performed during the computation.

\medskip

\noindent \underline{\textbf{Step degree 2}}~
Since $\operatorname{LM}_{\prec}(f_i)=x_1^2$ for $1\leq i \leq 4$, all six S-pairs among the input generators have sugar degree $2$.
In the direct Buchberger computation considered here, these pairs are processed at step degree $2$.
Put $g_1:=f_1$.
In the execution considered here, the pairs
$(f_1,f_2)$, $(f_1,f_3)$, and $(f_1,f_4)$ are processed first,
and their monic nonzero remainders are, respectively,
{\small
\begin{eqnarray*}
g_2&=&
\underline{x_1x_2}+41x_2^2+23x_1x_3+64x_2x_3+49x_3^2
 +16x_1+56x_2+57x_3,\\
g_3&=&
\underline{x_2^2}+14x_1x_3+43x_2x_3+22x_3^2+29x_3,\\
g_4&=&
\underline{x_1x_3}+3x_2x_3+33x_3^2+22x_1+5x_2+14x_3.
\end{eqnarray*}}

\noindent The remaining three S-pairs among $f_2,f_3,f_4$ reduce to zero after $g_2,g_3,g_4$ have been inserted.
Indeed, one has $S(f_2,f_3)=19g_2+21g_3$, $S(f_2,f_4)=25g_2+61g_3+30g_4$, and $S(f_3,f_4)=6g_2+40g_3+30g_4$.
Moreover, every new S-pair involving one of $g_2,g_3,g_4$ has sugar degree at least $3$.
Hence no further nonzero remainder is obtained at step degree $2$.
In particular,
\[
{\rm LM}_{\prec}(\mathcal G^{(<D)})
=
\{x_1^2,x_1x_2,x_2^2,x_1x_3\}
=
{\rm LM}_{\prec^h}(G_{\rm hom})_{<D}.
\]

\noindent \underline{\textbf{Step degree $D=3$}}~
At this point, we have $\operatorname{LM}_{\prec}(g_1)=x_1^2$, $\operatorname{LM}_{\prec}(g_2)=x_1x_2$, $\operatorname{LM}_{\prec}(g_3)=x_2^2$, and $\operatorname{LM}_{\prec}(g_4)=x_1x_3$.
The four nonzero remainders inserted at step degree $3$ are obtained, in this order, from the S-pairs $(g_2,g_4)$, $(g_1,g_4)$, $(g_2,g_3)$, and $(g_1,g_2)$.
The corresponding least common multiples are $x_1x_2x_3 \prec x_1^2x_3 \prec x_1x_2^2 \prec x_1^2x_2$, in accordance with the normal selection strategy used to break ties among S-pairs of the same sugar degree.
The resulting monic remainders are
{\small
\begin{eqnarray*}
g_5&=&
\underline{x_2x_3^2}+41x_3^3+5x_2x_3+35x_3^2
 +64x_1+42x_2+11x_3,\\
g_6&=&
\underline{x_3^3}+35x_3^2+37x_1+61x_2+24x_3,\\
g_7&=&
\underline{x_2x_3}+13x_3^2+3x_1+37x_2+72x_3,\\
g_8&=&
\underline{x_3^2}+72x_1+66x_2+70x_3.
\end{eqnarray*}}

\noindent
In a direct Buchberger computation, all other S-pairs of sugar degree $3$ reduce to zero.
Moreover, these zero reductions can be carried out without increasing the sugar degree, and every new S-pair involving one of $g_5,g_6,g_7,g_8$ has sugar degree at least $4$.
Thus $g_5,g_6,g_7,g_8$ are precisely the nonzero elements inserted at step degree $D=3$.

All four polynomials have sugar degree $D$.
Since the sugar degree cannot decrease along a reduction sequence, the reductions producing $g_5,g_6,g_7,g_8$ never increase the sugar degree.
Together with the preceding observation on the zero reductions, this shows that the modified computation of Section~\ref{subsec:overD} produces no additional intermediate polynomial at step degree $D$.
Hence $\mathcal G^{(D)}=\widehat{\mathcal G}^{(D)}=\{g_5,g_6,g_7,g_8\}$.

Among them, $g_7$ and $g_8$ have total degree $2$, and hence are
degree-fall elements.
After homogenizing them to sugar degree $D$, their leading monomials are $x_2x_3y$ and $x_3^2y$, respectively.
Consequently,
\[
{\rm LM}_{\prec}(\mathcal G^{(D)})
=
{\rm LM}_{\prec^h}(G_{\rm hom})_D\big|_{y=1}
=
\{x_2x_3^2,x_3^3,x_2x_3,x_3^2\}.
\]
Note that the degree-fall element $g_7$ is actually used as a reducer in the computation of $g_8$.
More precisely, the reduction of $S(g_1,g_2)$ contains the step $67x_2x_3+65x_3^2+58x_1+18x_2+15x_3$ $\xlongrightarrow[g_7]{}$ $
70x_3^2+3x_1+21x_2+9x_3$, whose monic form is $g_8$.
Since ${\rm sdeg}(g_7)=D$ and it is used without a monomial multiplier, this reduction preserves the sugar degree.
Thus the use of $g_7$ is consistent with the reduction convention of Definition~\ref{def:r(g1,g2)}.
This explicitly illustrates the reduction convention used in Lemma~\ref{lem:cor-atD}.

\medskip

\noindent \underline{\textbf{Step degree $4$}}~
Continuing the computation yields the following three nonzero remainders of sugar degree $4$:
\[
g_9=\underline{x_1}+61x_2+51x_3,\qquad g_{10}=\underline{x_2}+70x_3,\qquad g_{11}=\underline{x_3}.
\]
The remaining S-pairs of sugar degree $4$ reduce to zero, so these are the only nonzero remainders inserted at this step degree.
The degree-$4$ homogenizations of $g_9,g_{10},g_{11}$ inter-reduce to $\{x_1y^3,x_2y^3,x_3y^3\}$, which is the degree-$4$ part of $G_{\rm hom}$.
This additional correspondence at step degree $4$ is not asserted in general; the general synchronization result of Section~\ref{subsec:overD} is needed only through step degree $D$.

\medskip

After inter-reduction, the reduced Gr\"obner basis of $\langle F\rangle$ is therefore $G=\{x_1,x_2,x_3\}$.
This also agrees with the inter-reduction of $G_{\rm hom}|_{y=1}$.

\subsection{An example beyond cryptographic semi-regularity}\label{app:example-noncsr}

The preceding example illustrates the correspondence under cryptographic semi-regularity.
We conclude with an example showing that the injectivity condition governing this correspondence is strictly weaker than cryptographic semi-regularity.

Let $K$ be an arbitrary field, let $R=K[x_1,x_2,x_3]$ with the DRL ordering $x_1\succ x_2\succ x_3$, and set
\begin{eqnarray*}
f_1=x_1^3+x_2,\quad f_2=x_2^3,\quad f_3=x_3^3,\quad f_4=x_2^3+x_2x_3^2,\quad f_5=x_1x_2x_3
\end{eqnarray*}
with $d_1=\cdots=d_5=3$.
The highest-degree homogeneous parts are
\[
f_1^{\rm top}=x_1^3,\quad f_2^{\rm top}=x_2^3,\quad f_3^{\rm top}=x_3^3,\quad f_4^{\rm top}=x_2^3+x_2x_3^2,\quad f_5^{\rm top}=x_1x_2x_3.
\]
Then $\langle F^{\rm top}\rangle = \langle x_1^3,x_2^3,x_3^3,x_2x_3^2,x_1x_2x_3 \rangle$.
Indeed, the monomial $x_1^2x_3^2$ does not belong to this ideal, whereas every monomial of degree $5$ does.
Therefore we obtain $D=d_{\rm reg}(\langle F^{\rm top}\rangle)=5$.

On the other hand, $x_3f_2^{\rm top}+x_2f_3^{\rm top}-x_3f_4^{\rm top}=0$ gives a nonzero class in $H_1(K_\bullet(\bm F^{\rm top}))_4$.
Indeed, the corresponding Koszul cycle is $(0,x_3,x_2,-x_3,0)$, whereas the degree-$4$ part of the second Koszul module is zero, since all five generators have degree $3$.
Thus $H_1(K_\bullet(\bm F^{\rm top}))_4\neq0$.
Since $4<D$, the sequence $\bm F^{\rm top}=(f_1^{\rm top},f_2^{\rm top},f_3^{\rm top}, f_4^{\rm top}, f_5^{\rm top})$ is not cryptographic semi-regular.

The homogenizations are $f_1^h=x_1^3+x_2y^2$ and $f_i^h=f_i$ for $2\leq i\leq5$.
With respect to the DRL ordering $x_1\succ x_2\succ x_3\succ y$, the reduced Gr\"obner basis $G_{\rm hom}$ of $\langle F^h\rangle$ is
\[
\{
\underline{x_1^3}+x_2y^2,\,
\underline{x_2^3},\,
\underline{x_3^3},\,
\underline{x_1x_2x_3},\,
\underline{x_2x_3^2},\,
\underline{x_2^2x_3y^2}
\}.
\]
Thus the first minimal leading monomial divisible by $y$ is
$x_2^2x_3y^2$, of degree $5$.
By Lemma~\ref{lem:char-no-y}, we obtain $d_{\rm inj}=5=D$.
Therefore Condition~\ref{cond:inj} holds with $d_0=D$, although
$\bm F^{\rm top}$ is not cryptographic semi-regular.

For completeness, the corresponding degree fall can also be seen directly.
First, $S(f_2^h,f_4^h)=-x_2x_3^2$, so its monic remainder $g:=x_2x_3^2$ is inserted at step degree $3$.
The only new S-pairs of sugar degree $4$ are those of $g$ with $f_3$ and $f_5$, and both reduce to zero:
\[
S(g,f_3)=x_3g-x_2f_3=0,\qquad S(g,f_5)=x_1g-x_3f_5=0.
\]
Thus no further nonzero remainder is inserted before step degree $5$.
At step degree $5$, we have
\[
S(f_1^h,f_5^h) = x_2x_3f_1^h-x_1^2f_5^h = x_2^2x_3y^2.
\]
Correspondingly, in the affine computation,
\[
S(f_1,f_5)=x_2x_3f_1-x_1^2f_5=x_2^2x_3.
\]
The leading monomials of the current intermediate basis are
\begin{eqnarray*}
&& \mathrm{LM}_{\prec} (f_1)= x_1^3,\qquad \mathrm{LM}_{\prec} (f_2)=\mathrm{LM}_{\prec} (f_4)=x_2^3,\\
&& \mathrm{LM}_{\prec} (f_3)=x_3^3,\qquad \mathrm{LM}_{\prec} (f_5)=x_1x_2x_3,\qquad \mathrm{LM}_{\prec} (g)=x_2x_3^2,
\end{eqnarray*}
none of which divides $x_2^2x_3$.
Hence $S(f_1,f_5)$ is already reduced with respect to the current intermediate basis, and its nonzero remainder is $x_2^2x_3$ itself.
Since $ \deg(x_2^2x_3)=3<5$, a degree fall occurs at step degree $5=d_{\rm inj}$, in agreement with Corollary~\ref{cor:dinj-fall}.

\def\thesection{Appendix \Alph{section}}
\section{Generalized cryptographic semi-regularity}\label{app:gen}
\def\thesection{\Alph{section}}

As in the main body of this paper, let $R=K[x_1,\ldots,x_n]$ be the polynomial ring in $n$ variables over a field $K$, and let $\bm F=(f_1,\ldots,f_m)$ be a sequence of inhomogeneous polynomials of degrees $d_1,\ldots,d_m$.
Put $F=\{f_1,\ldots,f_m\}$ and $R'=R[y]$.

\subsection{Motivation and definition}

The following discussion was motivated by exploratory computations over finite fields.
We first consider the overdetermined case $m>n$; the case $m=n$ will be discussed separately in Remark~\ref{rem:mn} below.

We constructed systems having at least one prescribed affine zero as follows.
We chose a point $\bm a=(a_1,\ldots,a_n)\in K^n$, generated polynomials $g_1,\ldots,g_m\in R$ of degrees $d_1,\ldots,d_m$ by choosing their coefficients at random, and set $f_j=g_j-g_j(\bm a)$ for $1\leq j\leq m$.
Then $\bm a$ is a common zero of $F$, while $f_j^{\rm top}=g_j^{\rm top}$ for $1\leq j \leq m$.
We considered those instances for which $F$ has only finitely many zeros over the algebraic closure $\overline K$.

The computations were performed with Magma V2.25-3~\cite{MagmaBCP}, using in particular the built-in functions \texttt{Dimension}, \texttt{HilbertSeries}, and \texttt{GroebnerBasis}.
They were intended to explore the behavior described below and to motivate the following definition, rather than to provide a statistical study of its frequency.

This construction also gives a natural explanation for the two properties that repeatedly appeared in our computations.
Since the prescribed affine zero affects only the lower-degree parts of the input polynomials, the highest-degree homogeneous parts $\bm F^{\rm top}$ have the same distribution as those of the randomly generated polynomials $g_1,\ldots,g_m$.
Thus it is natural to compare their behavior with the usual cryptographic semi-regularity expected for random homogeneous systems.

On the other hand, the homogenized system $\bm F^h$ always has the projective zero $(a_1:\cdots:a_n:1)$.
Moreover, when $\bm F^{\rm top}$ is cryptographic semi-regular, it follows from $m>n$ that $R/\langle F^{\rm top}\rangle$ is Artinian, and hence $F^h$ has no projective zero with $y=0$.
Thus, for the instances considered above satisfying the property~{\rm (1)} below, $F^h$ has finitely many projective zeros.
It is then natural to ask whether the Hilbert function initially has the same behavior as that predicted for a homogeneous system with the given degree pattern and subsequently stabilizes at the number of projective zeros counted with multiplicity.

Put $D=d_{\rm reg}(\langle F^{\rm top}\rangle)$ and $D'=\widetilde d_{\rm reg}(\langle F^h\rangle)$.
The following two properties were repeatedly observed:
\begin{enumerate}
\item[(1)] $\bm{F}^{\rm top}$ is $D$-regular, i.e., cryptographic semi-regular;
\item[(2)] $\bm F^h$ is $D'$-regular.
\end{enumerate}
The property~{\rm (1)} is the usual cryptographic semi-regularity condition.
Over an infinite field, Fr\"oberg's conjecture~\cite{Froberg} predicts the corresponding Hilbert series for generic homogeneous polynomials of prescribed degrees.
We do not claim that the property~{\rm (2)} holds generically.

For these instances, we have $D'<\infty$, and the property~{\rm (2)} can be written as
\[
{\rm HS}_{R'/\langle F^h\rangle}(z)
=
\left(
\frac{\prod_{j=1}^m(1-z^{d_j})}
     {(1-z)^{n+1}}
\bmod z^{D'}
\right)
+
\sum_{d=D'}^\infty N_{F^h}z^d,
\]
where $N_{F^h}$ is the number of projective zeros of $F^h$, counted with multiplicity.
This motivates the following extension of cryptographic semi-regularity.

\begin{definition}
A homogeneous polynomial sequence $\bm H=(h_1,\ldots,h_m)$ is said to be \emph{generalized cryptographic semi-regular} if it is $\widetilde d_{\rm reg}(\langle\bm H\rangle)$-regular when $\widetilde d_{\rm reg}(\langle\bm H\rangle)<\infty$.
When $\widetilde d_{\rm reg}(\langle\bm H\rangle)=\infty$, this means that $\bm H$ is $d$-regular for every positive integer $d$.
\end{definition}

Every cryptographic semi-regular sequence is generalized cryptographic semi-regular.
In particular, the property {\rm (2)} above says that $\bm F^h$ is generalized cryptographic semi-regular.
We include the consequence below because it directly sharpens the degree bounds established in Section~\ref{sec:3}; a systematic study of generalized cryptographic semi-regularity and further consequences is undertaken in~\cite{KY2}.

\subsection{A refined bound under generalized cryptographic semi-regularity}

We next record a consequence of generalized cryptographic semi-regularity for the degree bounds considered in this paper.
For positive integers $n,d_1,\ldots,d_m$, recall the notation $D^{(n;d_1,\ldots,d_m)}$ from~\eqref{eq:newbound0}.

\begin{proposition}\label{prop:gen-csr-bound}
Assume that $m>n$, that $\bm F^{\rm top}$ is cryptographic semi-regular, and that $\bm F^h$ is generalized cryptographic semi-regular.
Put $D=d_{\rm reg}(\langle F^{\rm top}\rangle)$ and $D'=\widetilde d_{\rm reg}(\langle F^h\rangle)$.
Then $D=D^{(n;d_1,\ldots,d_m)}$ and
\begin{equation}\label{eq:Dnew}
\mathrm{max.GB.deg}_{\prec^h} (F^h) \leq \max \{ D,D'\} \leq  D^{(n+1;d_1,\ldots,d_m)}
\end{equation}
Consequently,
\[
\mathrm{sd}^{\rm mut}_{\prec}(F) \le \mathrm{sd}^{\rm mac}_{\prec}(F) \leq \mathrm{sd}^{\rm mac}_{\prec^h}(F^h)= \mathrm{max.GB.deg}_{\prec^h}(F^h) \le D^{(n+1;d_1,\ldots,d_m)}.
\]
\end{proposition}

\begin{proof}
Since $m>n$ and $\bm F^{\rm top}$ is cryptographic semi-regular, we have $D<\infty$ and $D=D^{(n;d_1,\ldots,d_m)}$.
Moreover, $R/\langle F^{\rm top}\rangle$ is Artinian, and the proof of Proposition~\ref{prop:new} shows that $D'<\infty$.
Since $\bm F^h$ is generalized cryptographic semi-regular, it is $D'$-regular.
By the definition of $D^{(n+1;d_1,\ldots,d_m)}$, the coefficient of degree $D^{(n+1;d_1,\ldots,d_m)}$ in $\frac{\prod_{j=1}^m(1-z^{d_j})}{(1-z)^{n+1}}$ is nonpositive.
If $D'>D^{(n+1;d_1,\ldots,d_m)}$, then $D'$-regularity would identify this coefficient with the corresponding value of the Hilbert function.
A negative value is impossible, while a zero value would force the Hilbert function to vanish in all subsequent degrees, and hence it would stabilize no later than $D^{(n+1;d_1,\ldots,d_m)}$, contrary to the definition of $D'$.
Hence $D'\leq D^{(n+1;d_1,\ldots,d_m)}$.
Moreover, we have $D=D^{(n;d_1,\ldots,d_m)}\le D^{(n+1;d_1,\ldots,d_m)}$, where the inequality follows directly from the definition.
Proposition~\ref{prop:new}~{\rm (2)} therefore yields \eqref{eq:Dnew}.
The assertions on the solving degrees now follow from the relations recalled in Section~\ref{subsec:knownbound}.
\end{proof}

Consequently, the arithmetic complexity over $K$ of computing a Gr\"obner basis of $\langle F^h\rangle$, and hence of $\langle F\rangle$ via homogenization, is bounded by
\begin{equation}\label{eq:new_comp_bound}
O\left(
m
\binom{
n+D^{(n+1;d_1,\ldots,d_m)}
}{
D^{(n+1;d_1,\ldots,d_m)}
}^{\omega}
\right),
\end{equation}
where $\omega$ denotes the exponent of matrix multiplication over $K$; see \cite[Theorem~1.1~(2) and Remark~3]{KY2}.

\begin{remark}\label{rem:mn}
When $m=n$, cryptographic semi-regularity of $\bm F^{\rm top}$ already implies that $\bm F^h$ is generalized cryptographic semi-regular and that $D'=D-1$; see \cite[Lemma~3.1.1]{KY2}.
Consequently, $\mathrm{max.GB.deg}_{\prec^h}(F^h)\leq D$.
Moreover, equality $\mathrm{max.GB.deg}_{\prec^h}(F^h)=D$ holds if the initial ideal $\langle\mathrm{LM}_{\prec^h}(\langle F^h\rangle)\rangle$ is weakly reverse lexicographic;
see Proposition~\ref{prop:new} and \cite[Theorem~1.2]{KY2}.
\end{remark}

The bound in~\eqref{eq:Dnew} can be substantially smaller than the first bound in Theorem~\ref{thm:maxGB} when the system is sufficiently overdetermined.
For example, in the quadratic case with $n=9$ and $m=18$, we have $D=D^{(9;2,\ldots,2)}=D^{(10;2,\ldots,2)}=4$, whereas the first bound in Theorem~\ref{thm:maxGB} is $11$.

The degree $D^{(n+1;d_1,\ldots,d_m)}$ also appears in analyses of algebraic attacks.
For example, \cite[Section~2.3]{FK24} uses $D^{(n+1;d_1,\ldots,d_m)}$ as a degree bound for the XL algorithm~\cite{XL}.
As discussed in \cite[Section~2.3, in particular Remark~1]{FK24}, this degree has also appeared in the cryptographic literature without an explicit regularity assumption on $\bm F^h$.
Generalized cryptographic semi-regularity gives a sufficient condition under which this degree bound, and consequently the arithmetic-complexity estimate~\eqref{eq:new_comp_bound}, is rigorously justified.

\end{document}